\documentclass[11pt,a4paper]{article}
\usepackage[pdftex]{graphicx}
\usepackage[cmex10]{amsmath}
\usepackage{amsfonts}
\usepackage{amssymb}
\usepackage{amsthm}
\usepackage{siunitx}

\usepackage{algorithm}
\usepackage{algpseudocode}
\usepackage{xcolor}
\usepackage[hidelinks=true]{hyperref}

\usepackage{xfrac}
\usepackage{enumerate}

\usepackage{subfigure}
\usepackage[font=footnotesize ]{caption}

\usepackage{multirow}
\usepackage{rotating}
\usepackage{makecell}
\usepackage{booktabs}
\usepackage{tabularx,ragged2e}
\usepackage{anysize}
\marginsize{2.5cm}{2.5cm}{2.5cm}{2.5cm}

\usepackage{setspace}
\usepackage{parskip}
\usepackage[pagewise]{lineno}

\usepackage{siunitx}

\usepackage{lipsum}

\usepackage{xcolor}

\definecolor{blue}{rgb}{0.0, 0.0, 0.0} 
\definecolor{red}{rgb}{0.0, 0.0, 0.0} 
\definecolor{teal}{rgb}{0.0, 0.0, 0.0} 

\usepackage{authblk}
\title{\textbf{{\color{blue}
Adaptive search space decomposition method for pre- and post- buckling analyses of space truss structures 
}}}
\author[1,4]{Varun~Ojha\thanks{
Corresponding Author: Varun~Ojha,~email: v.k.ojha@reading.ac.uk\\
Cite as: Ojha V, Pant{\`o}, and Nicosia G (2022) Adaptive search space decomposition method for pre- and post- buckling analyses of space truss structures \textit{Engineering Application of Artificial Intelligence}}
}
\author[2]{Bartolomeo Pant{\`o}}
\author[3]{Giuseppe Nicosia}

\affil[1]{University of Reading, Reading, UK}
\affil[2]{Durham University, Durham, UK}
\affil[3]{University of Catania, Catania, Italy}
\affil[4]{Newcastle University, Newcastle, UK}

\date{}
\begin{document}
\maketitle

\begin{abstract}
{\color{blue}
The paper proposes a novel adaptive search space decomposition method and a novel gradient-free optimization-based formulation for the pre- and post-buckling analyses of space truss structures. Space trusses are often employed in structural engineering to build large steel constructions, such as bridges and domes, whose structural response is characterized by large displacements. Therefore, these structures are vulnerable to progressive collapses due to local or global buckling effects, leading to sudden failures. The method proposed in this paper allows the analysis of the load-equilibrium path of truss structures to permanent and variable loading, including stable and unstable equilibrium stages and explicitly considering geometric nonlinearities. The goal of this work is to determine these equilibrium stages via optimization of the Lagrangian kinematic parameters of the system, determining the global equilibrium. However, this optimization problem is non-trivial due to the undefined parameter domain and the sensitivity and interaction among the Lagrangian parameters. Therefore, we propose formulating this problem as a nonlinear, multimodal, unconstrained, continuous optimization problem and develop a novel adaptive search space decomposition method, which progressively and adaptively re-defines the search domain (hypersphere) to evaluate the equilibrium of the system using a gradient-free optimization algorithm. We tackle three benchmark problems and evaluate a medium-sized test representing a real structural problem in this paper. The results are compared to those available in the literature regarding displacement-load curves and deformed configurations. The accuracy and robustness of the adopted methodology show a high potential of gradient-free algorithms in analyzing space truss structures.

}

\end{abstract}


\section{Introduction}
\label{sec:intro}
Space truss structures represent one of the most extensively used structural typologies in civil engineering to build non-ordinary steel structures such as bridges, large span arches, domes, and transmission towers. These structures generally show large displacements even under service loadings, and their ultimate response is often characterized by snap-through post buckling mechanisms in which the structure passes rapidly from an equilibrium state to a non-adjacent equilibrium configuration~\cite{hrinda2010snap}. Furthermore, some catastrophic events~\cite{martin2001another} and numerical studies~\cite{smith1984space,murtha1988alternate,blandford1996progressive} revealed that truss structures are prone to activate progressive collapses due to equilibrium instability, leading to brittle and sudden failure, which may lead to significant economic and human losses. Hence, complex nonlinear analyses considering the mechanical and geometrical nonlinearities that allow for a complete equilibrium path of the structure are needed to assess the structural robustness of truss structures against progressive collapses~\cite{gilbert2003layout}. 

In the last decades, many structural optimization procedures, based on linear programming~\cite{gilbert2003layout,tyas2006practical,horst2013global}, genetic algorithms~\cite{saka2007optimum,hasanccebi2009performance,huang2019engineering,lute2009computationally}, or iterative finite element procedures~\cite{hrinda2008optimization}, have been developed to research optimum design of large truss structures through exploring  optimum topology~\cite{yildiz2013comparison} or geometry of the system and optimum cross-sectional dimensions for the members. In these methods, the objective function has generally been the total weight of the structure, and the material tensile and compressive strengths have been considered as {\color{red}constraints}, eventually considering the elements' local buckling. However, these methods do not provide information on the nonlinear post-buckling response of the structure and its effective safety level. {\color{blue} Authors have also used neural networks and optimization algorithms for other structural engineering problems such as structural damage detection~\cite{tran2021efficient,khatir2021improved}, and 
simulation of fracture mechanics~\cite{khatir2019fast}, and monitoring of structural health~\cite{khatir2019structural}.}

Many authors proposed mathematical optimization algorithms based on mathematical programming to perform the limit analysis on structural systems such as frictional rigid-block assemblages~\cite{gilbert2006limit,ferris2001limit,baggio2000collapse} or elastoplastic Von-Mises steel structures~\cite{bisbos2007second,bisbos2005second}. These approaches provide valuable information about the failure mechanism and corresponding load multiplier. However, they are generally formulated under the hypothesis of small displacements and do not consider the incremental loading process. Thus, they cannot be applied for buckling analyses. 

On the other hand, incremental procedures within the finite element framework, based on the Arc-length method~\cite{crisfield1991nonlinear,crisfield1996nonlinear,memon2004arc} combined with iterative Newton-Raphson techniques~\cite{de2012nonlinear,bathe2006finite}, represent effective tools for performing nonlinear and buckling analysis considering the material and geometric nonlinearities. According to these strategies, an iterative stepwise linearization of the nonlinear structural behavior is considered. At each iteration, the tangent stiffness matrix and the geometrical stiffness matrix describe the mechanical and geometrical nonlinearities of the system, respectively. Thus, iterative Newton-Raphson procedures are effectively employed to evaluate the nonlinear response of 2D and 3D large structures, even in the presence of damage-plasticity constitutive laws~\cite{chen2007plasticity,macorini2011non}. However, they require a significant computational effort to update the stiffness matrices at each analysis iteration. Moreover, close to the critical points of the structural response, where the equilibrium configuration changes from stable to unstable (or vice-versa), numerical issues may appear, significantly increasing the number of iterations required to get the solution or leading to the divergence of the solution. 

A few learning algorithms for neural networks have been developed and applied for the nonlinear modeling of mechanical systems~\cite{memon2004arc} and for approximation of nonlinear behavior of structures~\cite{wang2015self}, covering the drawbacks of Newton-Raphson like procedures. These algorithms are based on quasi-Newton methods~\cite{li1995mechanical,geradin1981computational} that estimate the inverse Hessian of an objective function from the gradient to enable Newton-like optimization algorithms. Thus, not requiring the assembling of global stiffness matrices. Despite their potentialities, these approaches have not yet been exploited for the nonlinear assessment of large structural systems. 

{\color{blue}
The main goal of the paper is to propose a novel formulation of the nonlinear pre- and post-buckling analysis of space truss structure using gradient-free optimization algorithms. In this formulation, we subject a space truss structure (system) to large displacements and optimize the Lagrangian kinematic parameters (displacement) and the load multiplier of the system to guarantee the global equilibrium of the system. In our formulation, we define an objective function (the global equilibrium) in terms of global unbalance determined as the difference between the vectors of external and internal forces on the system. To the best of the authors' knowledge, this is the first work performing incremental (multistep) analysis of structures using gradient-free global optimization algorithms. 

The presented  optimization problem can be classified as a nonlinear, multimodal, unconstrained, continuous optimization problem. It is a challenging problem to solve as the search space for displacement and load multiplier variables has an undefined upper bound. Additionally, the search landscape poses significant challenges to existing continuous optimization algorithms. Hence, this problem could be considered as a testbench for optimization algorithms. The implementation of this research work is available on our GitHub page\footnote{\url{https://github.com/vojha-code/Hypershpere-Search}}. 

The displacement and load multiplier have their lower bound equal to zero. In this paper, we show experimental approaches considered to decompose and define search space in order to solve this optimization problem. Such experimental approaches allowed us to propose a novel search space decomposition method, \textit{adaptive search space decomposition}, which progressively and adaptively defines new hyperspheres (a bounded search space) for solving the optimization problem. The method is based on iteratively finding new centers for new hypersphere such that the centers follow the global equilibrium path a space truss structure. Thus allowing the study of the equilibrium stages and potential snap-through mechanisms of a space truss structure. 
}

In this work, three benchmark problems investigated in the literature through Newton-like approaches~\cite{crisfield1991nonlinear} are considered to assess the accuracy and performances of the selected set of optimization algorithms and proposed adaptive search space decomposition method. Each problem is solved using ad-hoc differential evolution algorithms (a gradient free optimization algorithms), comparing the results with those available in the literature regarding displacement load capacity curves and failure deformed configurations. Finally, a medium-size structure (a test problem on a 3D reticular beam), designed in the literature following the criteria of minimum weight, is solved to prove the applicability of the proposed procedure and to solve relevant problems in structural engineering.

{\color{blue}
The rest of the paper is organized as follows: First we introduce the structural model and its formulation as an optimization problem in Sections~\ref{sec:structural_model} and ~\ref{sec:opt_problem}. Then the search space decomposition methods are discussed in methodology Section~\ref{sec:method}. The adaptive search space decomposition method is described in Section~\ref{sec:hypersphere_algo}. The results on three} benchmark problems and a test problem are  presented and analyzed in Section~\ref{sec:results}. {\color{teal} Finally,} discussions and conclusions are presented in Section~\ref{sec:discussion} and Section~\ref{sec:conclusion}.

\section{\color{blue}Space truss structure problem}
\label{sec:the_structure}
\subsection{\color{blue} Kinematic and static of space truss structure}
\label{sec:structural_model}
\begin{figure}
    \centering
    \includegraphics{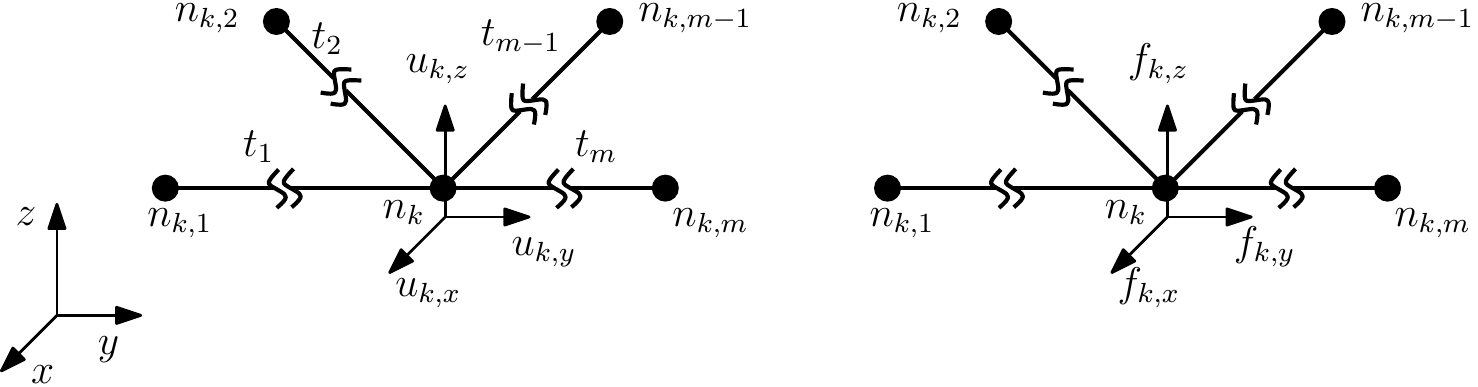}
    \caption{Independent displacements (\textit{Left}) and dual forces (\textit{Right}) of a generic node $n_k$ connecting $m$ trusses $t_1, t_2, \ldots, t_m$ with nodes $n_{k,1}, n_{k,2}, \ldots, n_{k,m}$. The node $n_k$ has its displacements $u_{k,x}, u_{k,y}, u_{k,z}$ and corresponding forces $f_{k,x}, f_{k,y}, f_{k,z}$  along $x, y$ and $z$ dimensions of the space.}
    \label{fig:ind_disp}
\end{figure}

\begin{figure}
    \centering
    \includegraphics[width=.55\textwidth]{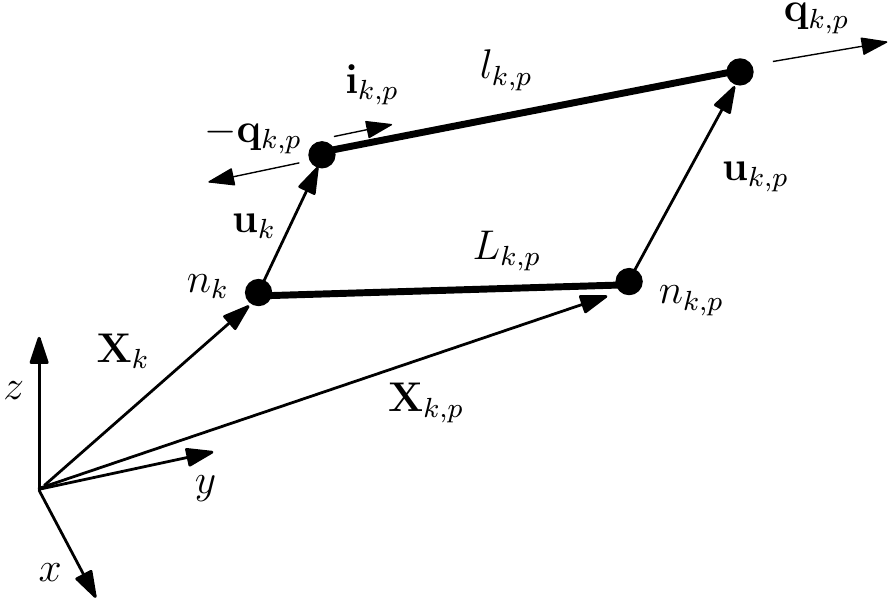}
    \caption{Kinematic and internal forces of a generic $p$-th truss connecting node $n_k$ with node $n_{k,p}$. The truss has an initial length $L_{k,p}$ and its nodess initial coordinates are $\textbf{X}_k$ and $\textbf{X}_{k,p}$. The axial nodal force of the truss is $\textbf{q}_p$, its  current length and orientation are $l_{k,p}$ and $\textbf{i}_{k,p}$ and nodal displacements are $\textbf{u}_k$ and $\textbf{u}_{k,p}$.}
    \label{fig:kinematics}
\end{figure}

The system kiematics is bescribed by model adopt a Lagrangian description by considering large displacements and small strain hypotheses. The degrees of freedom are assumed {\color{blue} to be} coincident with the 3$N$ absolute displacements of the $N$ free nodes of the truss system, referred to as the global reference system on x-y-z axes (Figure~\ref{fig:ind_disp}). The displacements of the generic node $n_k$ with $k = 1, \ldots, N$, are collected in a vector $\textbf{u}_k = [u_{k,x}, u_{k,y}, u_{k,z}]$ while the corresponding dual nodal forces are collected in a vector $\textbf{f}_k = \textbf{f}_{0k} + \lambda \textbf{f}_k [f_{0k,x}, f_{0k,y}, f_{0k,z}] + \lambda [f_{k,x}, f_{k,y}, f_{k,z}]$ where $\textbf{f}_{0k}$ represents the permanent loads and $\textbf{f}_k$ represents variable loads, which are amplified by the load multiplier $\lambda$ (Figure~\ref{fig:ind_disp}). Each node can be connected with an arbitrary number $(m)$ of nodes $(n_{k,1}, \ldots, n_{k,m})$ by as many trusses $(t_1, \ldots, t_m)$ as shown in Figure~\ref{fig:ind_disp}. The generic $p$-th truss $(p = 1, \ldots, m)$ connecting the node $n_k$ with the node $n_{k,p}$ is represented in Figure~\ref{fig:kinematics}, where the vector $\textbf{X}_k$ and $\textbf{X}_{k,p}$ represent the initial coordinates of the truss end nodes, $\textbf{u}_{k,p}$ the displacement vector of the node $n_{k,p}$ and $\textbf{q}_{k,p}$ the {\color{red} nodal} axial forces {\color{blue}acting on the truss}. 

According to the hypothesis of small {\color{red}deformations} \cite{crisfield1991nonlinear}, the {\color{teal}truss deformation} ($\epsilon_p$) is evaluated considering the engineering strain {\color{teal} definition} as: 
\begin{equation}
    \delta_{k,p} = \frac{l_{k,p}}{L_{k,p}} - 1
\end{equation}
where  $L_{k,p} = ||\textbf{X}_{k,p} - \textbf{X}_k ||$ and $l_{k,p} = || \textbf{X}_{k,p} + \textbf{u}_{k,p} - \textbf{X}_k - \textbf{u}_k ||$ are the initial and current length of the truss, respectively and {\color{teal}$|| \cdot ||$ indicates the Euclidean norm}.

Given $k_p$ the axial stiffness of the truss and $\textbf{i}_{k,p}$ the unitary vector identifying its current orientation, the internal force can be written as follows:
\begin{equation}
    \color{teal}
    \textbf{q}_{k,p}  = k_p \delta_{k,p} \textbf{i}_{k,p} = k_p \left( \frac{|| \textbf{X}_{k,p}- \textbf{X}_k||}{|| \textbf{X}_{k,p} + \textbf{u}_{k,p}-\textbf{X}_k - \textbf{u}_k||} - 1)\right)\textbf{i}_{k,p},
\end{equation}
where
\begin{equation}
    \textbf{i}_{k,p}=\frac{(\textbf{X}_{k,p}+\textbf{u}_{k,p}-\textbf{X}_k-\textbf{u}_k)}{\color{teal}l_{k,p}}.
\end{equation}
The equilibrium between the external and internal forces at the generic node $n_k$ can be expressed as follows: 
\begin{equation}
    \label{eq:equilibrium_generic_node}
    \textbf{f}_{0k} + \lambda \textbf{f}_k + \sum_{p=1}^m \textbf{q}_{k,p} = \textbf{0}~\quad  \text{for each } k = 1, \ldots, N.
\end{equation}
{\color{blue} Note that this work does not consider material nonlinearities (i.e., a linear elastic constitutive law characterizes structural elements). Therefore, truss stiffnesses assume constant values during the analysis. 
}

\subsection{\color{blue}The global equilibrium as an optimization problem}
\label{sec:opt_problem}
{\color{teal} In this Section, the global equilibrium of the truss system is mathematically characterized and formulated as an optimization problem whose free variables are represented by the Lagrangian parameters of the system, collected in the vector} $\textbf{u} = [\textbf{u}_1, \ldots, \textbf{u}_N ]$  and the load multiplier $\lambda$. The objective function $(\mathcal{F})$ {\color{teal} that the optimization algorithms aim to minimize} is expressed in terms of global unbalance as follows: 
\begin{equation}
    \label{eq:main_objective}
    \mathcal{F}(\textbf{u},\lambda) =\frac{\sqrt{\textbf{R}_1^2 +\textbf{R}_2^2+ \cdots \textbf{R}_N^2}}{\sqrt{\textbf{f}_1^2+ \textbf{f}_2^2+ \cdots +\textbf{f}_N^2}} = 0,
\end{equation}
where $\textbf{R}_k, k=1,\ldots,N$  is the unbalance vector of the $k$-th node, which can be easily obtained from Equation~\eqref{eq:equilibrium_generic_node}. 
{A value 0 is the global minima of Equation~\ref{eq:main_objective} indicates an equilibrium state of the structure.}

\begin{figure}
    \centering
    \includegraphics[width=.6\textwidth]{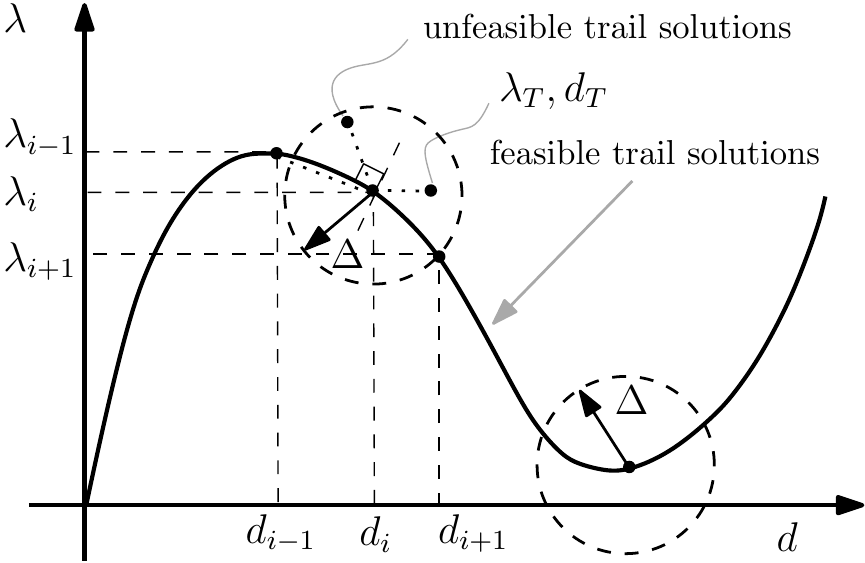}
    \caption{Multistep procedure. Variables $d$, $\lambda$, and $ \Delta $ respectively are control point, load multiplier, and discrete steps along the control point dimension for producing an equilibrium path (curve in black). Any {\color{teal} optimal point (a solution that satisfy  Equation~\ref{eq:main_objective}; otherwise sub-optimal solution)} on this curve and within a radius $ \Delta $ that satisfies constraints in Equations~\eqref{eq:const1} and \eqref{eq:const2} is a feasible solution {\color{blue} in terms of its characteristics being within a tolerance proximity to the equilibrium path formed by a vertical hypothesis load. Solutions not within this tolerance proximity defined by Equations~\eqref{eq:const1} and \eqref{eq:const2} are unfeasible or belong to a path concerning non-vertical hypothesis load.}}
    \label{fig:Multi_Step}
\end{figure}

The incremental multistep procedure shown in Figure~\ref{fig:Multi_Step}, {\color{teal}known in the literature as arc-length method~\cite{crisfield1991nonlinear}, can be adapted} to draw the equilibrium path of the structure in the $d-\lambda$ space, where $d$ is a chosen control point, which can coincide with the physical displacement of a structural node or a function of a number of displacements. For instance, the control point can be assumed coincident with the Euclidean norm of the vector $\textbf{u} (d = ||\textbf{u}||)$. {\color{blue} This multistep arc-length procedure is framed within the equilibrium optimization problem, by considering the equality constraint reported in Equations~\ref{eq:const1} and  \ref{eq:const2}, where} $d_i$ and $\lambda_i$ {\color{red} are the control point and load multiplier at the generic $i$-th equilibrated configuration belonging to the equilibrium path,} and $d_T - \lambda_T$ the corresponding trial values:
\begin{equation}
    \label{eq:const1}
    \sqrt{((d_T - d_i )^2 + (\lambda_T - \lambda_i)^2 )} = \Delta,
\end{equation}
where $\Delta$ is a parameter of the procedure governing the discretization of the equilibrium path; furthermore, if $i >  1$, a second constrain is considered as:
\begin{equation}
    \label{eq:const2}
    (d_T - d_i ) (d_{i-1}-d_i) + (\lambda_T-\lambda_i ) (\lambda_{i-1}-\lambda_i ) \le 0,
\end{equation}
where, $d_{i-1}$ and $\lambda_{i-1}$ are the values of the control point and load multiplier at step $i-1$.

\subsection{\color{blue}Challenges of the global equilibrium optimization}
\label{sec:opt_problem_challenges}

{\color{teal}The challenging issue with optimization problem defined in Section~\ref{sec:opt_problem} is the lack of a properly defined search domain for the free variables displacement and load multiplier, to be pursued using gradient-free optimization algorithms. The search domain for the free variables is only the initial elastic prediction of the structural response and simple equilibrium considerations.
}

{\color{red}The unavailability of a straightforward search domain definition for the space truss structure  makes it an ill-posed problem. However, the gradient-free optimization algorithms allow experimenting with an intuitive guess of the search domain. Hence, we formulate this space truss structure buckling analysis problem as an optimization problem.} 

The geometric nonlinearities, high physical interaction and sensitivity among displacement variables, and numerous possible imbalance states for a particular load multiplier make this problem further challenging. {\color{blue} This can be seen later in results Section~\ref{sec:results} mentioned as sub-optimal and local minima solutions.} 

We, therefore, classify the optimization of Equation~\eqref{eq:main_objective} as a nonlinear, multimodal, unconstrained, continuous, minimization problem whose search landscape is complex containing many local and global minima, i.e., its search landscape is rugged.

In order to solve this optimization problem, first, an intuitive guess was made based on an initial elastic prediction of the structural response and simple equilibrium considerations. However, each displacement variable has its own domain, and fixing one min-max domain range to all variables complicates the search landscape, leading to an impossible optimization (demonstrated later in Section~\ref{sec:results}). {\color{teal} Therefore, a number of strategies for search space decomposition discussed in Section~\ref{sec:method}  are used in this research to solve this problem as efficiently as possible.
}

{\color{blue} In this paper, we show experimental approaches considered to decompose and define search space in order to solve this optimization problem. Furthermore, such experimental approaches allowed us to propose a novel adaptive search space decomposition method for multistep analyses, denoted as the hypersphere method, representing an efficient and robust transposition of the arc-length method into the optimization framework.  }

\section{\color{blue}Methodology}
\label{sec:method}
The procedure adopted in this study to solve the optimization problem defined in Section \ref{sec:the_structure}, is based on gradient-free algorithms and search space partition methods characterized by different levels of complexity. First, an \textit{informed decomposition} of search space was based on the equilibrium-path curve (also referred to as baseline equilibrium path), available for some selected benchmarks, already investigated in the literature by adopting Newton-Raphson iterative procedures combined with arc-length methods~\cite{hrinda2010snap,crisfield1991nonlinear,de2012nonlinear}. Second, a more sophisticated \textit{adaptive search space partition} technique is used. 

\subsection{\color{blue}Informed decomposition of search space}
\label{sec:search_domain}
{\color{teal} In this study, we consider two alternative} options for decomposing the search space {\color{teal}for this problem. First, we} partition {\color{red} the control point displacement target} into several discrete steps as the literature provided the information on the control point. 

Second, in addition to partitioning the control point, {\color{teal} we incrementally partition the space for all free} displacement variables. Compared to the first option, the second option was more effective (as discussed later in Section~\ref{sec:results}) in finding solutions across the {\color{blue} baseline} equilibrium path shown in the literature~\cite{crisfield1991nonlinear,crisfield1996nonlinear}. The performance of each optimization algorithm is discussed later in the results reported in Section~\ref{sec:results}.  

{\color{teal} These search domain partition strategies highlighted the strengths and weaknesses of the optimization algorithms in solving this challenging class of optimization problems.} The goal of the space truss structures optimization was to find solutions across the known equilibrium path of the benchmark problems as dense as possible. Therefore, \textit{denseness}, \textit{coverage of the path}, the \textit{number of optimal solutions}, and \textit{speed of the convergence}, among others, were the main criteria for evaluating the {\color{teal} effectiveness of our methodology} in solving this nonlinear, unconstrained, continuous and  multimodal, optimization problem.

\subsection{\color{blue}Adaptive decomposition of search space} 
\label{sec:hypersphere_algo}
\begin{figure}
    \centering
    \includegraphics{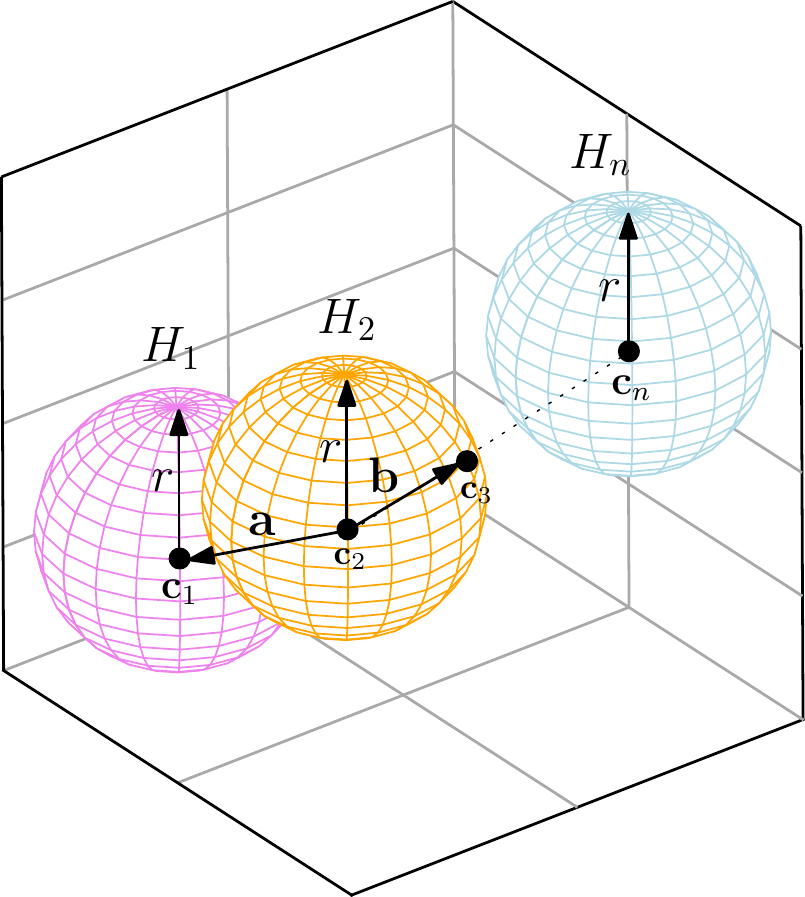}
    \caption{Adaptive search space decomposition method (i.e., hypersphere search algorithm). An initial hypersphere $H_1$ is first constructed based on an input radius $r$ and the center $\textbf{c}_1$. Hyperspheres $H_2, \ldots, H_n$ are discovered by the algorithm automatically by taking a boundary point (i.e., a solution within a hypersphere that has the largest distance from the previous hypersphere's center). For example, center $\textbf{c}_2$ has a distance $ \textbf{a} $ from center $ \textbf{c}_1 $ and center $\textbf{c}_3$ has a distance $ \textbf{b} $ from center $ \textbf{c}_2 $.
    	\label{fig:hypersphere}}
\end{figure}

{\color{blue} Instead of decomposing the problem search space in a pure discrete form mentioned in Section~\ref{sec:search_domain}, a method for \textit{adaptive decomposition of search space} is proposed in this Section. In this method, we automatically search hyperspheres that define suitable search space for the multi-step pre-and post buckling analysis of space truss structures.

This is an iterative procedure where the method starts with an initial seed hypersphere. The initial seed hypersphere is defined by a very small range of the free displacement and load multiplier variables. This is to simulate a minute displacement of the structure, which is observed through the control point (node at which vertical load is applied) dimension.

For this seed hypersphere, the method at first evaluates a few trail solutions (e.g. five solutions). Then from all optimal solutions, it searches for a new hypersphere.} This process continues until the max displacement of the control point of a problem is reached (see Figure~\ref{fig:hypersphere}).

{\color{blue} The proposed adaptive search space decomposition method, called  \textit{hypersphere search algorithm}, finds the most suitable hyperspheres to support gradient-free optimization algorithms to find the global equilibrium effectively. The method comprises the following steps:}
\begin{enumerate}
    \item[]  \textbf{Step 1.} The first step is to initialize a hypersphere {\color{teal}$H_k$ with center in $\textbf{c}_k$ for $k = 1$} and radius $r$. {\color{teal}The search of the centers is controlled by a 2-dimension vector $[d,~\lambda]$, where $d$} is the initial displacement control point (generally zero) and $\lambda$ is the corresponding load multiplier that guarantees the global equilibrium solution. A gradient-free optimization algorithm then evaluates {\color{blue}one or more trial solutions within this hypersphere.} 
	
	\item[]  \textbf{Step 2.} The second center {\color{teal}$\textbf{c}_k$ for $k = 2$} is the farthest point in the previous hypersphere {\color{blue}$H_{k-1}$} from its center {\color{blue}$\textbf{c}_{k-1}$} among all trial solutions in Step 1 that {\color{teal} are optimal (a value 10$^{-5}$)}. {\color{blue}In other words, the second center $\textbf{c}_k$ is set to the solution that has max value of $d$ and is optimal.} Then the hypersphere {\color{teal}$H_k$ for $k = 2$} is formed around the center {\color{teal}$\textbf{c}_k$} with radius $r$. Similar to Step 1, one or more  trial solutions within this new hypersphere $H_k$ are evaluated for obtaining the hypersphere {\color{teal}$H_{k+1}$}.
	
	\item[]  \textbf{Step 3.} The {\color{teal}next center $\textbf{c}_k$ for $k=3$} is the farthest point from the center {\color{teal}$\textbf{c}_{k-1}$} among all trial solution points $\textbf{c}_i$ obtained in the hypersphere {\color{teal}$H_{k-1}$} and is the point that satisfies the following condition:
	\begin{equation}
        \textbf{b}_i \cdot \textbf{a}  < 0  \quad \text{and} \quad ||\textbf{b}_i||  > ||\textbf{b}_j||,\quad i \ne j \text{ for all } i \text{ and } j \text{ in } H_{k-1},  	    
	\end{equation}
	where $i$ and $j$ are indexes of trial solutions in the hypersphere {\color{teal}$H_{k-1}$, vector $\textbf{a} = \textbf{c}_{k-2} - \textbf{c}_{k-1}$, and $\textbf{b}_i =  \textbf{c}_i - \textbf{c}_{k-1}$. Once the center $\textbf{c}_k$ is obtained, a hypersphere $H_k$} is formed with a radius $r$.  This step is shown in Figure~\ref{fig:hypersphere}.
	
    \item[]  \textbf{Step 4.} {\color{blue}Stop if $k$ is reached a max number of hypersphere trails or a center reached a predefined max displacement $d_{\text{max}}$ is reached. Else set $k = k + 1$ go to Step 3 to find other hyperspheres.}
\end{enumerate}

The user-defined hyper-parameter of this algorithm is the radius $r$, the number of trial solutions needed to be produced in a hypersphere, and the stopping criteria. While the radius $r$ depends on the user set value, the number of trial solution evaluation depend on the experimentation. The stopping criteria can be set to either an arbitrary sufficient number of solutions that may approximately cover equilibrium path obtained by the arc-length method, or to a maximum displacement value of the control point dimension obtained by arc-length method (or an intuitive guess {\color{teal}that} could be made about the buckling of the structure).

~\\
\subsection{Experiment setup}
\label{sec:exp_setup}
{\color{blue}
As described in Section \ref{sec:method}, the optimization is performed by adopting gradient-free algorithms and search of the solution was supported by space partitioning methods. In this Section, the setups of these approaches are described. This methodology is validated on three benchmarks already investigated in the literature by means of iterative Newton-Raphson procedures arc-length methods and tested on a medium-sized test problem. The three benchmarks consisted in 3D shallow truss structures, for which the baseline is available from the literature, enabling the setting of the min-max displacement and load multiplier value. However, these min-max values were an intuitive guess of the range, i.e., the maximum displacement that the control point displacement can take and the range of variability of external load. Therefore, the domain range was initially set to the min and max values of these benchmark problems in the literature. However, optimization algorithms could not find solutions for the first two benchmark problems across the equilibrium path shown in the literature, as discussed later in results Section 4. Hence, a strategy to decompose the search space was adopted. Finally, the results obtained for the three considered benchmarks were compared with the baseline from the literature.   

}
\subsubsection{\color{blue} Gradient-free algorithms for the equilibrium path analysis}
\label{sec:opt_algo}
{\color{teal} This section briefly explains each optimization algorithm used for the optimization.} The first set of algorithms chosen were single solution optimization algorithms DIRECT and Simulated Annealing (SA). Both algorithms work on a single solution, and they can provide a quick understanding of the problem to be optimized. DIRECT is a pattern search algorithm of a class of Lipschitz optimization methods that relies on a dividing rectangle principle, hence it is called DIRECT~\cite{jones1993lipschitzian}. SA, a heuristic-based gradient-free algorithm~\cite{kirkpatrick1983optimization}, probabilistically accepts solutions while exploring the neighborhood of a search point. The probability of acceptance depends on an energy factor called temperature that goes to zero from a certain initial value in a controlled manner by a factor reducing it in each iteration. This is analogous to the annealing process in metallurgy, hence the name Simulated Annealing. 

The second class of algorithms were swarm intelligence-inspired population-based optimization algorithms: Artificial Bee Colony Optimization (ABC), Ant Colony Optimization (ACO), and Particle Swarm Optimization (PSO). These algorithms emulate the foraging behavior of the swarm. For example, ABC follows the foraging behavior of honeybees~\cite{karaboga2008performance}, while ACO follows ants~\cite{socha2008ant}, and PSO follows the behavior of a flock of birds or a school of fish~\cite{kennedy1995particle,clerc1999swarm}. Each of these algorithms has its own framework and convergence properties. For example, ABC maintains three types of honeybees: employed, onlooker, and scout bees. ACO uses Gaussian distribution for its pheromone and solution matrix update. PSO uses an inertia factor~\cite{kennedy1995particle} (PSO-Std) or constriction factor~\cite{clerc1999swarm} (PSO-Const) to update the velocity and position of particles.  However, common features are updating a solution vector for a defined objective function that governs the solution's quality that needs an update in every generation of the optimization.

The final class of algorithms and the main focus of this research is the differential evolution (DE) algorithm from the family of evolutionary algorithms, which uses mutation, crossover, and selection operations on the solution vectors~\cite{storn1997differential}. There are several versions of the DE algorithm~\cite{opara2019differential}. After an initial performance evaluation (in terms of speed of convergence and diversity of solutions in multiple runs) of a set of 12 versions of DE strategies available in Scipy~\cite{scipyde}, this paper selects the two most promising DE versions: DE/rand/1/bin and DE/best/2/bin. These two versions vary in how they select a base vector [randomly (rand) or best] and the number of differences of vectors (one difference vector or two differences of vectors) for the mutation operation in DE. However, both do a binomial (bin) crossover. 


\subsubsection{\color{blue} Setup of hyperparameters of the algorithms}
\label{sec:algo_param_setting}
Initial trials on benchmark problems were performed to fixate hyperparameter values of the algorithms. For DIRECT, except for termination criteria, there was no other hyperparameter to set. {\color{teal} Moreover, DIRECT is a single solution based optimization algorithm. Similar to DIRECT, SA is a single solution optimization algorithm.} However, SA has its other hyperparameters temperature and initial temperature reduction rate was set to 0.1 and 0.99.

The swarm-based algorithms population size was set to 50, where other hyperparameters specific to each algorithm were as follows: {\color{teal} the hyperparameters of ABC} was abandonment limit and acceleration coefficient upper bound and they were set to 0.6 and 1; {\color{teal} the hyperparameters of ACO} sample size, selection pressure, and deviation-distance ratio, respectively set to 40, 0.5, and 1; {\color{teal} hyperparameters of PSO-Std} inertia weight, inertia weight damping ratio, personal learning coefficient, global learning coefficient respectively set to 1, 0.99, 1.5, and 2.0; and {\color{teal} hyperparameters of PSO-Const} personal learning coefficient and global learning coefficient were equal to \textit{inertia weight}$ \times \phi $, where inertia weight $ = 2/(\phi-2+\sqrt{\phi^2-4\times\phi})$ for $\phi $ and inertia weight damping set to 2.05 and 1. Finally, the  hyperparameters of DE algorithm scaling factor bound and crossover rate were set to [0.2, 0.9] and 0.9. Two versions of mutation ``DE/rand/1/bin'' and ``DE/best/2/bin were used. 

The initial radius for the adaptive search space decomposition method (see Section~\ref{sec:hypersphere_algo}) was set to 5 and termination condition was set to max displacement in the baseline of the problems. The rest of the experiment setting {\color{teal} about benchmark problems are described in results in Section~\ref{sec:results} as when a they are introduced.}

\section{Results and analysis}
\label{sec:results}

\subsection{Benchmark 1: Eight-member shallow truss structure}
\label{sec:benchmark_1}
\begin{figure}
        \centering
        \includegraphics[width=0.4\linewidth]{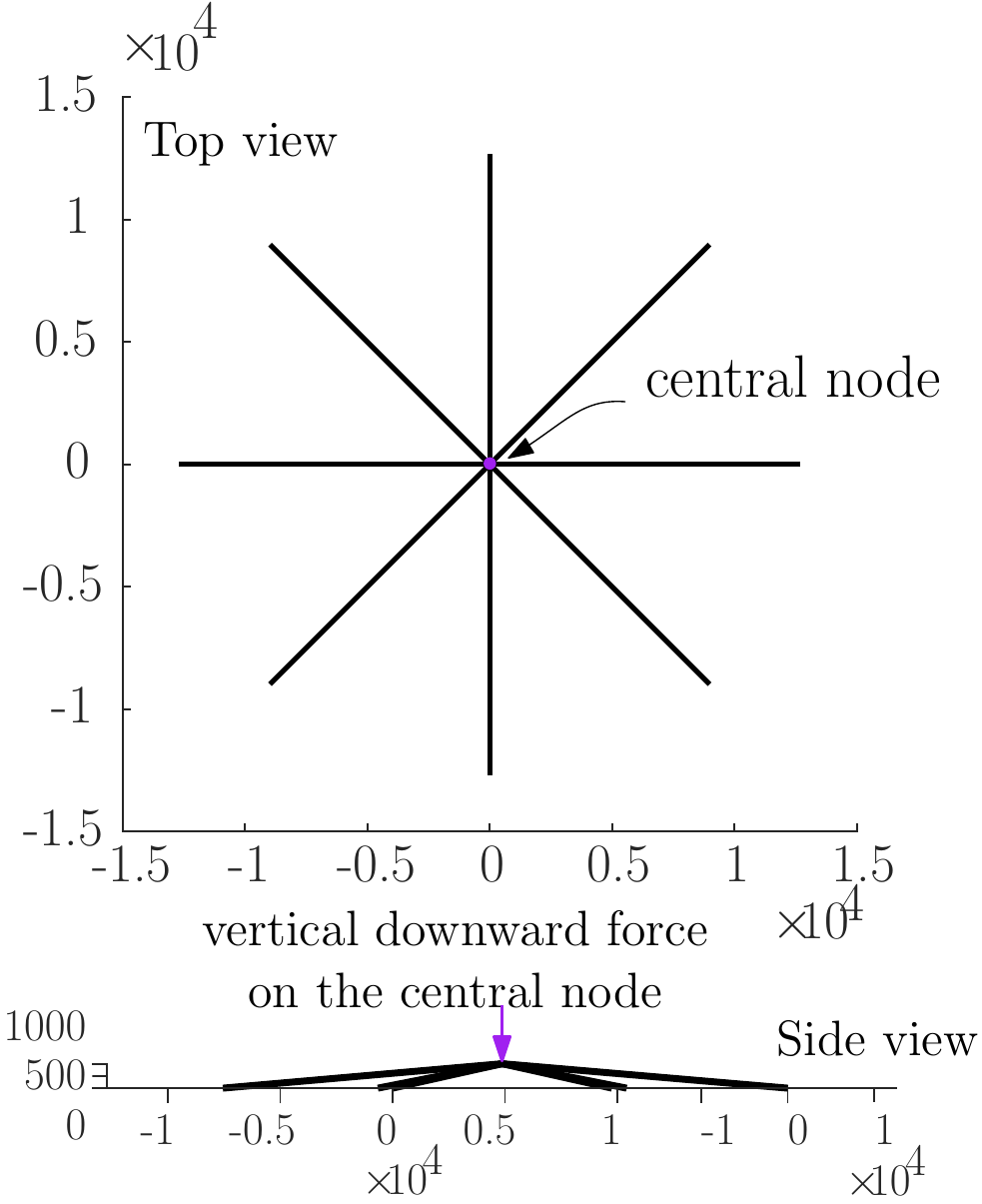}
        \caption{Eight-Member Shallow Truss Structure: Top view and Side view (dimensions in mm). The central node of the structure is indicated by a lighter color dot in Top view and a light color arrow pointer in Side view. The arrow on Side view indicates a vertical downward external force of \SI{4.45}{\kN} applied on the central node.
            \label{fig:8_member}}
\end{figure}

The first benchmark problem was a shallow truss structure comprising eight members fully restrained at the base and connected to each other by a central node (Figure~\ref{fig:8_member}). Each truss has a horizontal length of \SI{12700}{\mm} and a vertical rise of \SI{1000}{\mm}. All the trusses have the same cross-section area, equal to \SI{6450}{\mm^2}, and Young's modulus of \SI{70}{\MPa}. The structure is loaded with a concentrated vertical force of \SI{4.45}{\kN} applied at the central node.

\begin{figure}
        \centering
        \includegraphics[width=0.98\linewidth]{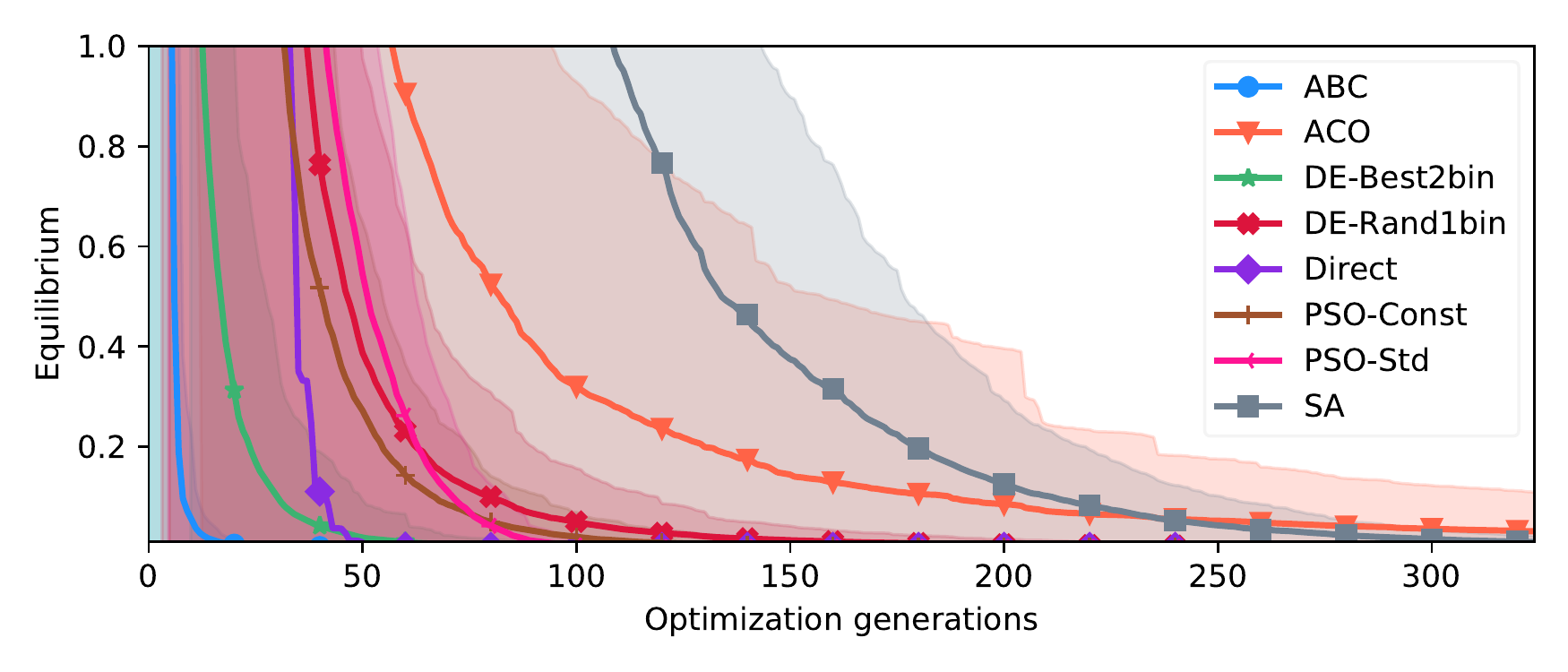}
        \caption{Convergence profile of algorithms averaged over 1000 trail solutions, all with optimum equilibrium value \num{1e-5}. Algorithms ABC, ACO, DE/best/2/bin, DE/rand/1/bin, Direct, PSO-Const, PSO-Std, and SA, respectively are in blue, orange, green, red, purple, brown, magenta, and gray; and they are respectively marked with symbols \(\bullet , \blacktriangledown , \bigstar ,  \times, \blacklozenge, |, <,\) and \( \blacksquare \). The shaded width of color around a line represents the standard deviation of the convergence.  \label{fig:eight_convergence}}
\end{figure}

The search space for the eight-member truss structure has three degrees of freedom (DoF) corresponding to the three translations of the central point ($N = 1$). Therefore, four variables were subjected to optimization, including the control point of the displacement variable and one load multiplier variable. The domain range for the experiments for the displacement was [\SI{0}{\mm}, \SI{3000}{\mm}], and for the load multiplier, it was set to [\num{-0.2}, 1]. All select algorithms were applied to optimize the eight-member truss structure on the full domain for initial trials. {\color{red}Initially, \num{1000} solutions were evaluated}. That is, the algorithms were run for \num{1000} instances. All algorithms for all \num{1000} trials converged to a global optimum (i.e., the precision of the order of \num{1e-5}. The convergence speed and variance of the algorithms are shown in Figure~\ref{fig:eight_convergence}.

\begin{figure}
        \centering
        \includegraphics[width=0.9\linewidth]{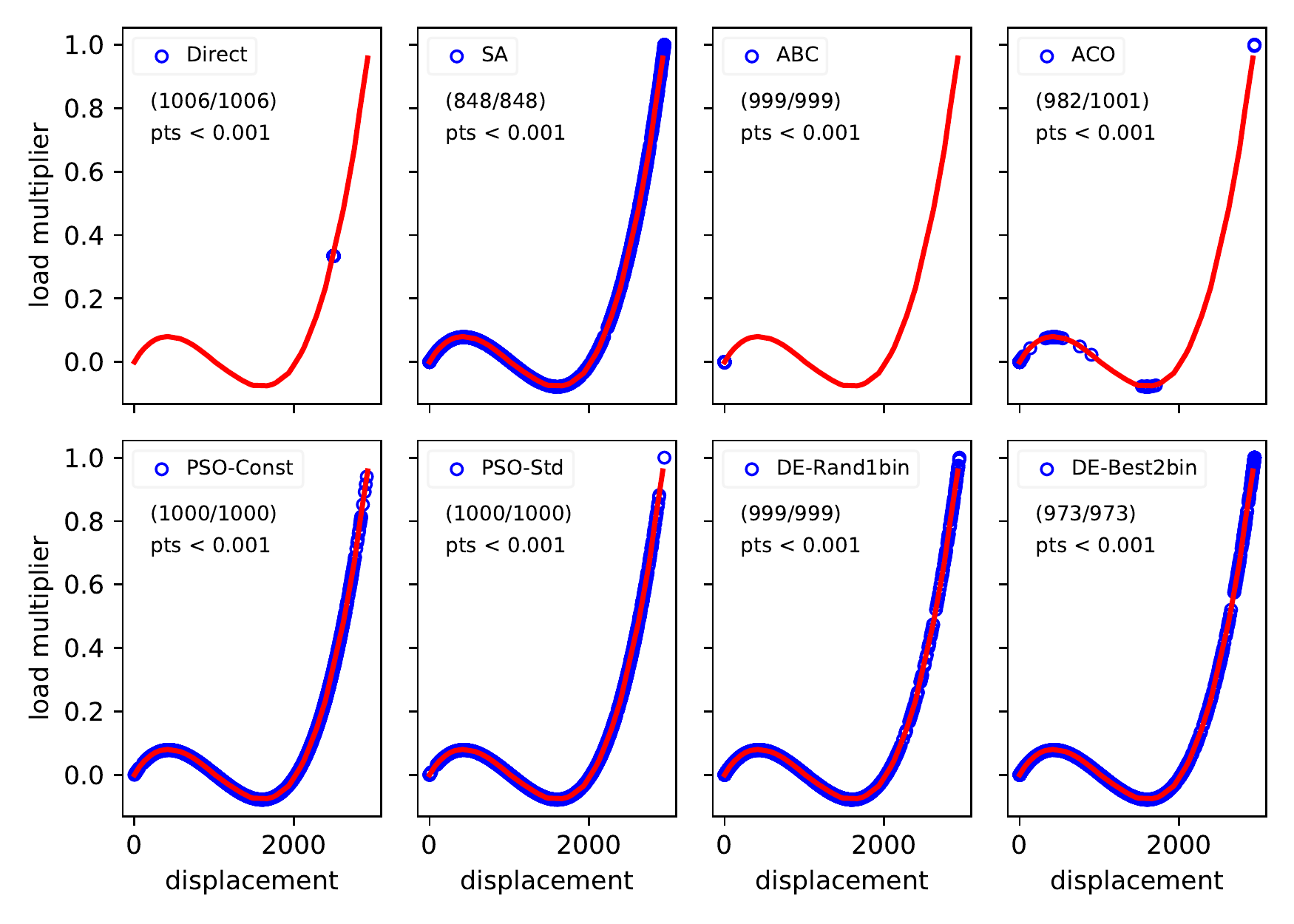}
        \caption{Eight-Member truss structure optimized solutions generated by searching with a domain setting of \([0, 1]\) for load multiplier and \([\SI{0}{\mm},  \SI{3000}{\mm}]\) for displacement variables (DoF), including the control point. Each plot is an x-y plane of control point displacement (x-axis) against load multiplier (y-axis). The solutions are shown by blue circles, and the red line represents the equilibrium path obtained by the standard Arc-length method mentioned in Section~\ref{sec:opt_problem}. From top left to bottom right, the algorithms are Direct (\textit{top left}), SA, ABC, ACO, PSO-Std, PSO-Const, DE/rand/1/bin, and DE/best/2/bin (\textit{bottom-right}). For each algorithm, the percentage of successful solutions (pts $ < 0.001$) generated are indicated in the plots, e.g., DIRECT produced 100\% successful points. i.e., \num{1006} solutions have a precision of 0.001 out of a total of \num{1006} trials. Similarly, SA produced 848/848 solutions with a precision of 0.001. 
        \label{fig:eight_algo_point}}
\end{figure}

\begin{figure}
        \centering
             \includegraphics[width=0.4\linewidth]{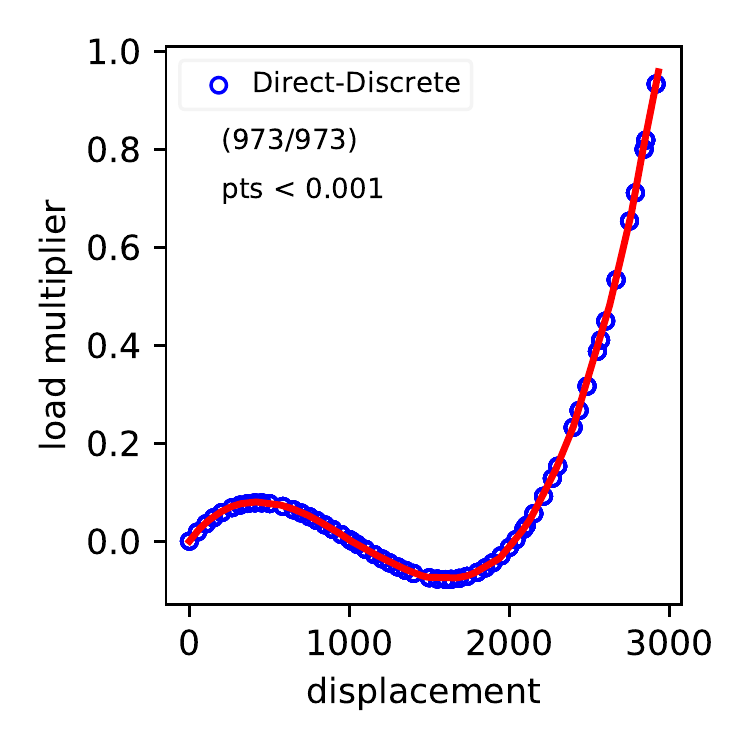} 
            \includegraphics[width=0.4\linewidth]{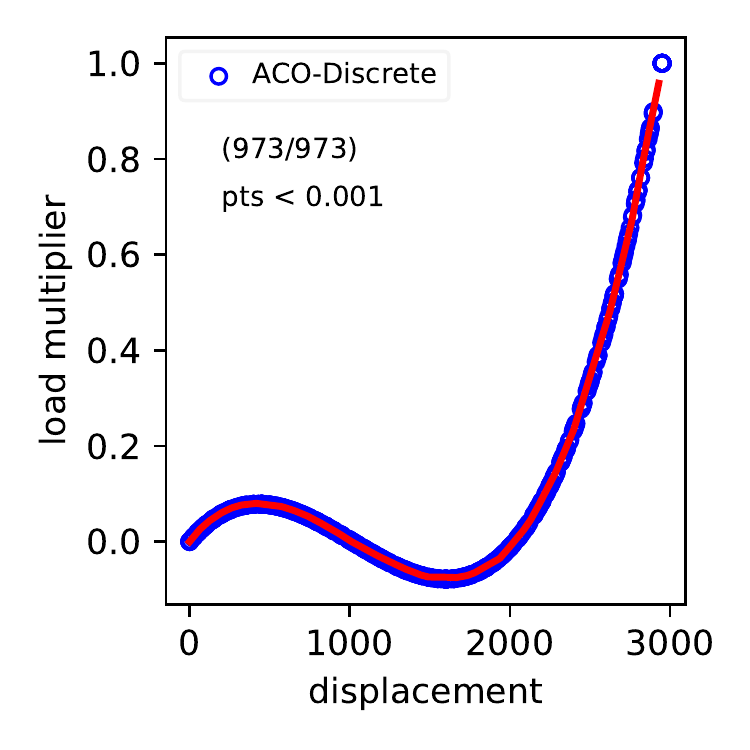}
        \caption{Solutions of Direct (\textit{Left}) and ACO (\textit{Right}) generated by decomposing domain of only control point (variable related to central node), which produced discrete sets like \([0, 250], [250, 500], \cdots, [2750, 3000]\). For each set, approx 100 trail solutions were {\color{red}evaluated}.
        \label{fig:eight_disc_point}}
\end{figure}

Figure~\ref{fig:eight_convergence} shows that ABC converged the fastest, followed by DE/best/2/bin, DIRECT, PSO-Std, PSO-Constriction (PSO-Const), DE/rand/1/bin, ACO, and SA. However, this performance is significant only if the algorithms were able to obtain solutions across the equilibrium path---Figure~\ref{fig:eight_algo_point} plots solutions obtained by the algorithms on a control-point and load multiplier axis to verify solutions. ABC, DIRECT, and ACO were not able to discover solutions across the equilibrium path. However, the performances of DIRECT and ACO were better than the performance of ABC. ABC failed to find any solution on the path except for a value of \num{0.0}. The other algorithms, SA, PSO-Const, PSO Std, DE/best/2/bin, DE/rand/1/bin, successfully covered the whole equilibrium path in \num{1000} trials. The continuous line reported in the graphs {\color{red} of Figure~\ref{fig:eight_algo_point} and \ref{fig:eight_disc_point}} represents the analytical solution from the literature~\cite{hrinda2010snap,crisfield1991nonlinear,de2012nonlinear}. 

The properties and the framework of the algorithms played a role in such characteristics of the solutions. For example, the DIRECT algorithm divides the rectangle and starts from the initial range, and partitions are based on the first solution obtained for the objective function for a given range. Moreover, since it has a deterministic approach toward optimization, it seemingly finds a solution at the center of the value of the first variable (load multiplier). Hence, in all \num{1000} trials, it found the exact same solution every time.

The other best-performing algorithms SA, PSO-Const, PSO Std, DE/best/2/bin, and DE/rand/1/bin, use uniform distribution to initialize the solutions. This has enabled them to cover the full range of the domain uniformly. On the other hand, ACO's new solutions are sampled from a Gaussian distribution, leading to the solutions following the central tendency (see Figure~\ref{fig:eight_algo_point}).

Since DIRECT and ACO were unable to cover the whole of the equilibrium path, the discrete control-point search space trials were performed to evaluate if these algorithms have the potential to find all other solutions for this problem. Both algorithms were able to cover the equilibrium path (see Figure~\ref{fig:eight_disc_point}). However, ACO could densely cover the equilibrium path, whereas DIRECT shows sparsity for the higher displacement and load multiplier values. These initial results helped select fewer best-performing algorithms to optimize the second benchmark's problem of the space truss structures.   

\subsection{Benchmark 2: Sixteen-member shallow truss structure}
\label{sec:benchmark_2}
\begin{figure}[H]
        \centering
             \includegraphics[width=0.3\linewidth]{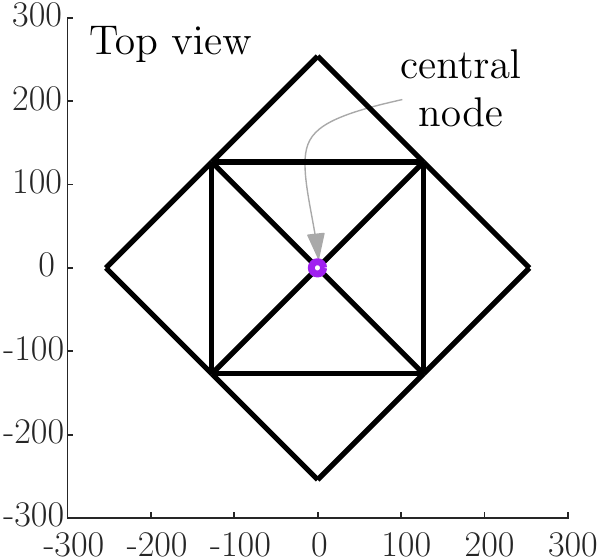} \quad
            \includegraphics[width=0.6\linewidth]{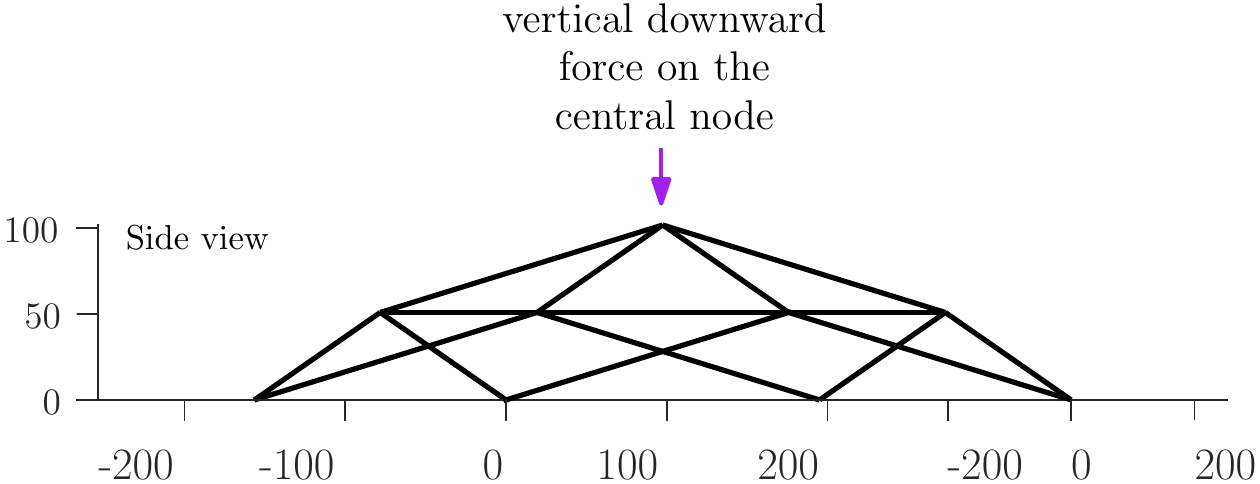}
        \caption{Sixteen-Member Shallow Truss Structure. \textit{Left:} Top View. \textit{Right:} Side View. Th central node of the structure is indicated by a lighter color dot in Top view and a light color arrow pointer in Side view. The arrow on Side view indicates a vertical downward external force of \SI{4450}{\kN} applied on the central node.
            \label{fig:3d_star_view}}
\end{figure}

The second benchmark problem was the Sixteen-Member Shallow Truss Structure (Figure~\ref{fig:3d_star_view}). The structure has four {\color{red}constraint} nodes and five free nodes ($N = 5$). Therefore, it is characterized by 15 DoF. The global dimensions in the plane are \SI{254}{\mm} in each direction, and the global vertical rise is \SI{100}{\mm}. All the trusses have the same cross-section area of \SI{645}{\mm^2} and Young's modulus of \SI{68950}{\MPa}. The structure is loaded with a concentrated vertical force of \SI{4450}{\kN} applied at the central node of the structure in Figure~\ref{fig:3d_star_view}.   

\subsubsection{Results of \color{blue} the informed search space decomposition methods}
\label{sec:star_search_domain_anal}
The search space for the sixteen-member truss structure comprises 16 variables to be optimized, including the 15 free displacements of the structure and the load multiplier variable. The domain range of the displacement variables was [\SI{0}{\mm}, \SI{250}{\mm}], and the load multiplier domain was set to [\num{-0.4}, 0.85].

\begin{figure}
        \centering
             \includegraphics[width=0.5\linewidth]{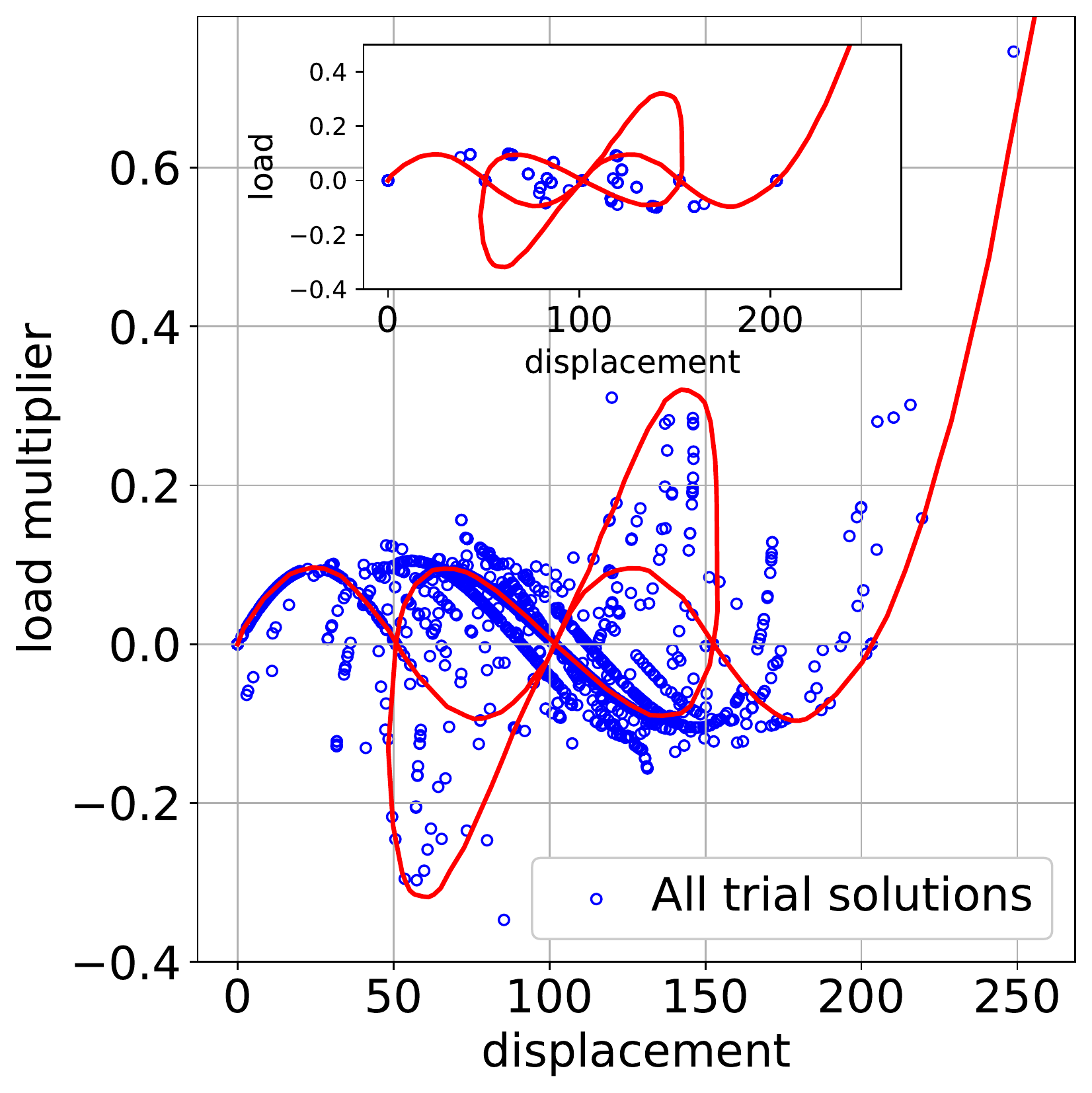}\hfill
            \includegraphics[width=0.5\linewidth]{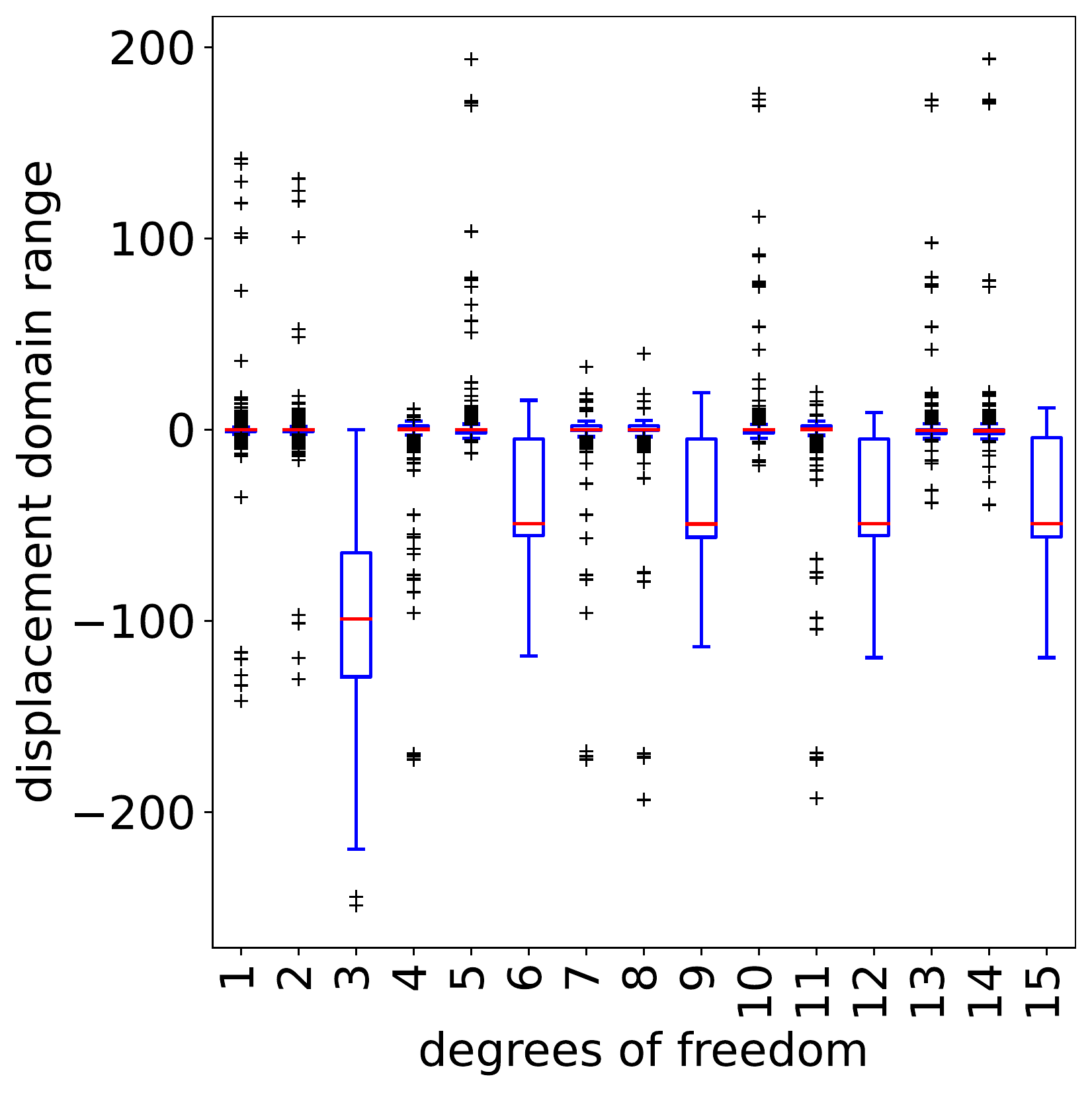}
        \caption{Sixteen Members initial trial solutions. \textit{Left:} initial trials for all discretized domain (inner-plot shows solution obtained for the experiment on full non-discretized domain [\SI{0}{\mm}, \SI{250}{\mm}]). \textit{Right:} Displacement variables  (DoF) domain analysis and approximation. \label{fig:3d_star_full_disc_domain_anal}}
\end{figure}

\textbf{Full domain analysis.} After setting an intuitive domain of the variables of the sixteen-member truss structure problem, the best performing algorithm of the first benchmark, i.e.,  DE/rand/1/bin was first applied to solve this sixteen-member truss structure. However, contrary to the results of the first benchmark, only a few solutions were obtained in all \num{1000} trials, and a tiny number of solutions were obtained on the equilibrium path (see the inner plot in Figure~\ref{fig:3d_star_full_disc_domain_anal}). Therefore, a decomposed domain analysis procedure (see Section~\ref{sec:search_domain}) was launched to find a suitable domain for respective variables. Again, the continuous red line in Figure~\ref{fig:3d_star_full_disc_domain_anal} indicates the analytical result from the literature, obtained under the hypothesis of vertical displacement of the central point.

\textbf{Domain decomposition analysis.} The decomposed domain analysis procedure partitioned the control point domain and the domain of other displacement variables. Based on the knowledge from the literature, the discretization of the control point dimension was straightforward. However, the discretization of the other variables could not be performed individually. Hence, an intuitive incremental guess was applied.

The grid space in Figure~\ref{fig:3d_star_full_disc_domain_anal}(left) shows the procedure adopted to analyze suitable domains for these variables. For each grid space, 50 trial solutions were {\color{red}evaluated}. All trial solutions that provided global-optimum value are plotted on an x-y plane of control point displacement and load multiplier in Figure~\ref{fig:3d_star_full_disc_domain_anal}(left). The DE/rand/1/bin provided several solutions within the range [\SI{0}{\mm}, \SI{225}{\mm}] of the control point in good agreement with the analytical solution available in the literature. Other solutions, far from the analytical curve, are provided by this optimization algorithm as well. However, these solutions correspond to alternative equilibrium paths characterized by a non-vertical displacement of the central node.

Figure~\ref{fig:3d_star_full_disc_domain_anal}(right) shows the displacement value distribution of each variable. The plot in Figure~\ref{fig:3d_star_full_disc_domain_anal}(right) provides the min-max range of each displacement variable for which global optimum solutions could be obtained. The knowledge of this min-max range (i.e., displacement domain knowledge) was then be used by all optimization algorithms.

\subsubsection{\color{blue}Results of the optimization on a learned search space}
\begin{figure}
        \centering
        \includegraphics[width=0.5\linewidth]{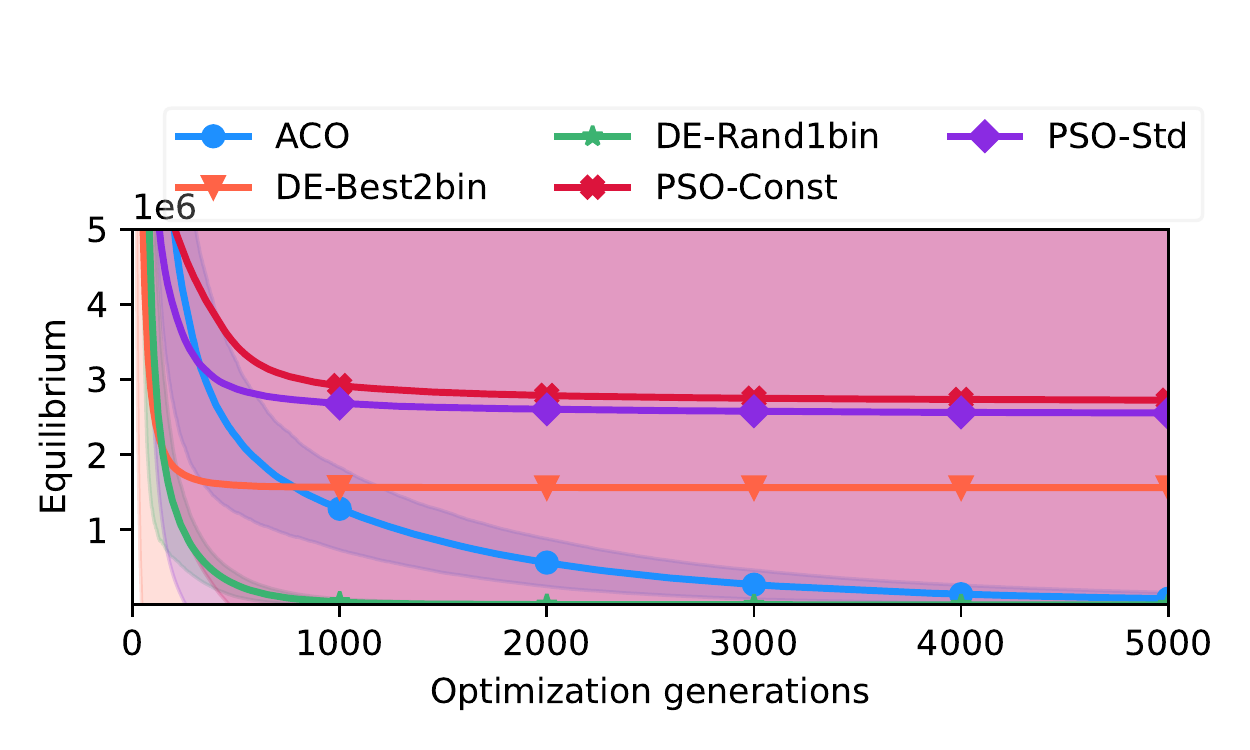}\hfill
        \includegraphics[width=0.5\linewidth]{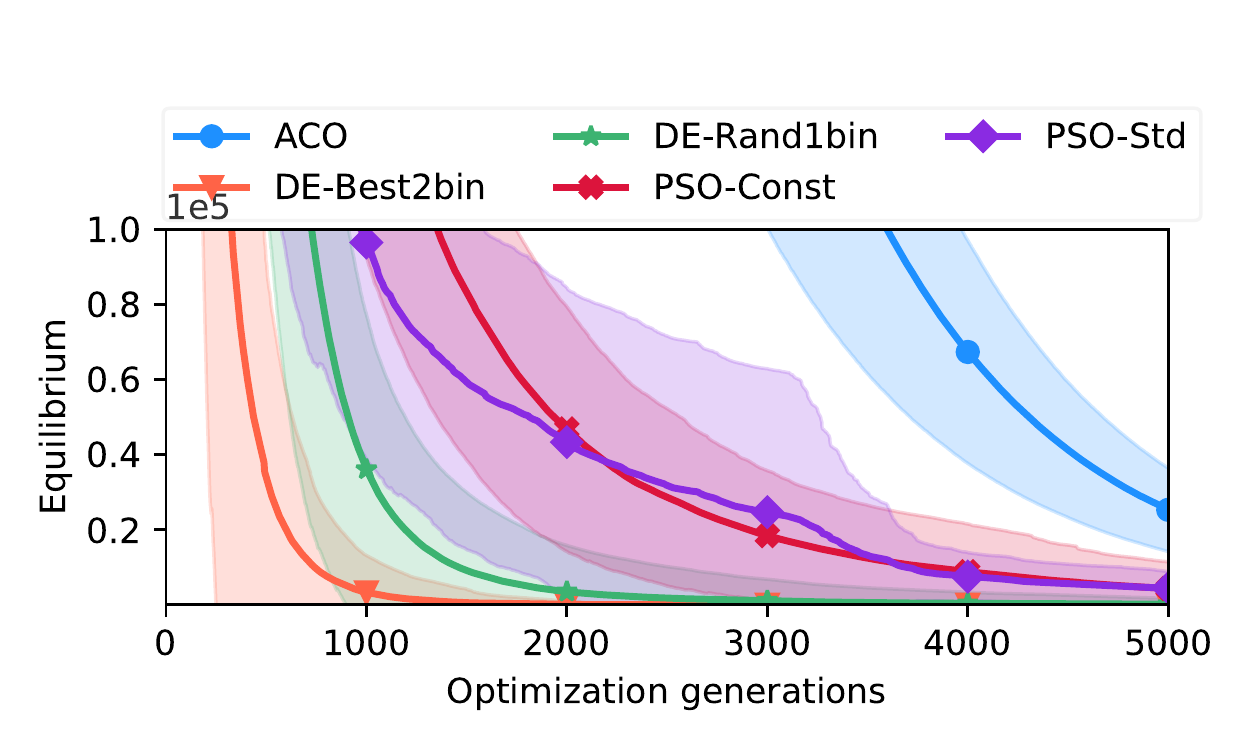}
        \caption{Optimization of Sixteen-Member shallow truss structure using five algorithms: (\textit{Left}) Convergence profile of algorithms averaged over solutions with equilibrium value of 500 and (\textit{Right}) a closer snapshot of convergence profile with non-converging solutions filtered out. Algorithms ACO, DE/best/2/bin, DE/rand/1/bin, PSO-Const, and PSO-Std are in blue, orange, green, red, and purple; and they are marked with symbols \(\bullet , \blacktriangledown , \bigstar ,  \times, \blacklozenge \). The shaded width of color around a line represents the standard deviation of the convergence. The x-axis indicate number of generations and y-axis indicate objective as per Equation~\eqref{eq:main_objective}. \label{fig:3d_star_convergence}}
\end{figure}

On the obtained knowledge of the domain shown in Figure~\ref{fig:3d_star_full_disc_domain_anal}(right), ACO, DE/best/2/bin, DE/rand/1/bin, PSO-Const, and PSO-Std algorithms were applied (with the domain setting of each displacement variable) for stopping criteria of solution accuracy \num{1e-5} and maximum iterations of \num{50000}. Figure~\ref{fig:3d_star_convergence}(left) shows each algorithm's average (of all solutions, including optimal and suboptimal) convergence profile. The average convergence profile rejects PSO-Const and PSO-Std algorithms. However, this only means that PSO obtained many poor (suboptimal) solutions (or even non-converging solutions) that affected the average convergence profile. A similar result was obtained for DE/best/2/bin, where poor solutions affected its average convergence profile. ACO and DE/rand/1/bin show more robustness in this scenario. They show the ability to escape local minima in most cases, which is observed from these two algorithms' smooth average convergence profile in Figure~\ref{fig:3d_star_convergence}(left).

Figure~\ref{fig:3d_star_convergence}(right) offers a microscopic view of the convergence of these algorithms. In Figure~\ref{fig:3d_star_convergence}(right), all solutions converged to a suboptimal equilibrium function [Equation~\eqref{eq:main_objective}] value of 500, selected after a few analyses of suboptimal points to see how best the algorithm's profile could be studied. The convergence profile is the average of all converging solutions to this suboptimal value, and all non-converging (all solutions stuck to local minima or did not converge at all) were filtered out. In this plot, DE/best/2/bin is the fastest converging algorithm, with all its select solutions being converged to global minima, i.e., \num{1e-5} in less than 5000 iterations.

Similarly, all solutions of the DE/rand/1/bin algorithm converged to global minima and this algorithm was the second-fastest. PSO-Const and PSO-Std were the subsequent fastest converging algorithms, where PSO-Const showed smoother convergence than the PSO-Std. However, the average convergence profiles of both reached near a suboptimal equilibrium value of 500 in 5000 iterations. Finally, the ACO algorithm's convergence was the last among this set of algorithms, and it is the slowest to reach the suboptimal value of 500.

\begin{figure}
        \centering
        \includegraphics[width=0.98\linewidth]{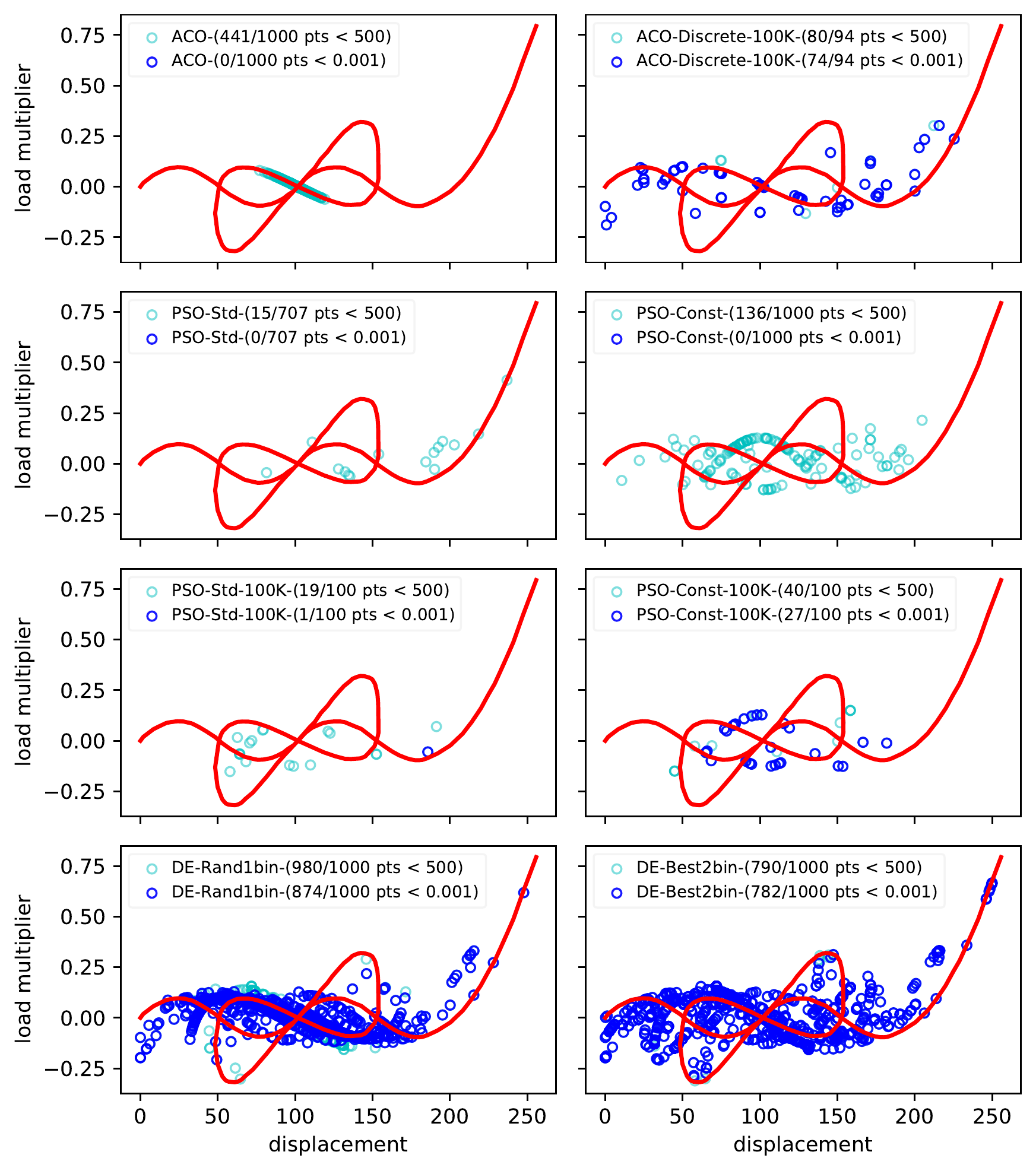}
        \caption{Quality of solutions (in blue solutions with a precision of 0.001 as their equilibrium value and solution in cyan are with precision 500) produced by optimization algorithms over \num{10000} generations and \num{100000} generations for ACO-Discrete and PSO-100K versions. Arc-length method based equilibrium path is shown in red. From top-left to bottom-right algorithms are ACO, ACO (with a discrete domain of control points over \num{100000} generations), PSO-Std, PSO-Const, PSO-Std with \num{100000} generations, PSO-Const with \num{100000} generations, DE/rand/1/bin, and DE/best/2/bin. \label{fig:3d_star_algo_points}}
\end{figure}

The convergence profiles analysis compared these algorithms in terms of their ability to escape local minima, speed of convergence, and the ability to find global optimum solutions. However, the accuracy of solutions of this class of space truss structure problem is also about how the solution compared to the analytical solution from literature like arc length~\cite{hrinda2010snap,crisfield1991nonlinear,de2012nonlinear} for achieving the equilibrium path characterized by a vertical displacement of the central node. Therefore, the filtered solution obtained by ACO, DE/best/2/bin, DE/rand/1/bin, PSO-Const, and PSO-Std in Figure~\ref{fig:3d_star_convergence}(right) were plotted on an x-y plane of control point displacement and the load multiplier in Figure~\ref{fig:3d_star_algo_points}.  

\textbf{Solutions accuracy analysis against the analytical curve.} Figure~\ref{fig:3d_star_algo_points} plots  suboptimal solutions  in cyan and optimal solutions in blue. As shown in Figure~\ref{fig:3d_star_algo_points}(row 1), the ACO algorithm, although accurate in finding solutions on the equilibrium path, was only able to find suboptimal solutions. However, the discrete domain analysis along the control point displacement dimension was performed for generating 100 more trial solutions of ACO with \num{100000} iterations as its stopping criteria. This analysis improved the results, and there were several optimal solutions on the equilibrium path; however, they were mainly scattered across the path. 

Similar to the ACO algorithm, PSO-Std and PSO-Const were unable to find any near-optimal solution in \num{50000} iterations. Additionally, PSO-Std and PSO-Const solutions were scattered (mostly belong to alternative equilibrium paths characterized by a non-vertical displacement), and rarely solutions were on the equilibrium path characterized by vertical displacement (see cyan points in row 2 of Figure~\ref{fig:3d_star_algo_points}). Since PSO version solutions were spread across the full domain of the problem, unlike ACO, it was not required to discretize the domain for further analysis. Instead, an analysis involving generating additional 100 new trial solutions for a \num{100000} iteration as stopping criteria was performed. Increasing iteration was for testing PSO robustness was evident from the fact that the solutions spread across, and the solutions were closer to suboptimal value than ACO (see Figure~\ref{fig:3d_star_convergence}(right)). However, even in this analysis, PSO-Std could not find any near-optimal solutions, and although PSO-Const found many near-optimal solutions, they were mainly inaccurate when compared with the required equilibrium path characterized by vertical displacement.

In contrast to ACO and PSO algorithms, DE versions were able to find a high percentage of accurate, near-optimal, and optimal solutions (Figure~\ref{fig:3d_star_algo_points}, last row). Interestingly, the DE version DE/rand/1/bin was more accurate than DE/best/2/bin. A possible explanation is that premature convergence to local minima is more frequent in DE/best/2/bin than in DE/rand/1/bin. This is because solutions of DE/best/2/bin start following the local best solution, and if the best solution finds a trajectory leading to local minima on the hypersurface of the solution search space, then the DE/best/2/bin will not converge to global minima. This fact is evident in Figure~\ref{fig:3d_star_convergence}(left), where DE/best/2/bin average convergence profile is poorly affected by solutions that are stuck into local minima.

Since both DE/rand/1/bin and DE/best/2/bin were able to provide a high percentage of near-optimal and accurate solutions, it was not necessary to do either a discrete or an extensive iteration analysis. However, many solutions were not on the equilibrium path characterized by vertical displacement. Rather, many solutions belong to the equilibrium path of varied characterizations. This is a disadvantage as it is only possible to pick a solution and analyze the equilibrium, but the space truss structure's equilibrium path (or a buckling direction) could not be observed clearly for the hypothesis's vertical downward force. Hence, an equilibrium path identification analysis was performed.

\begin{figure}
        \centering
        \includegraphics[width=0.98\linewidth]{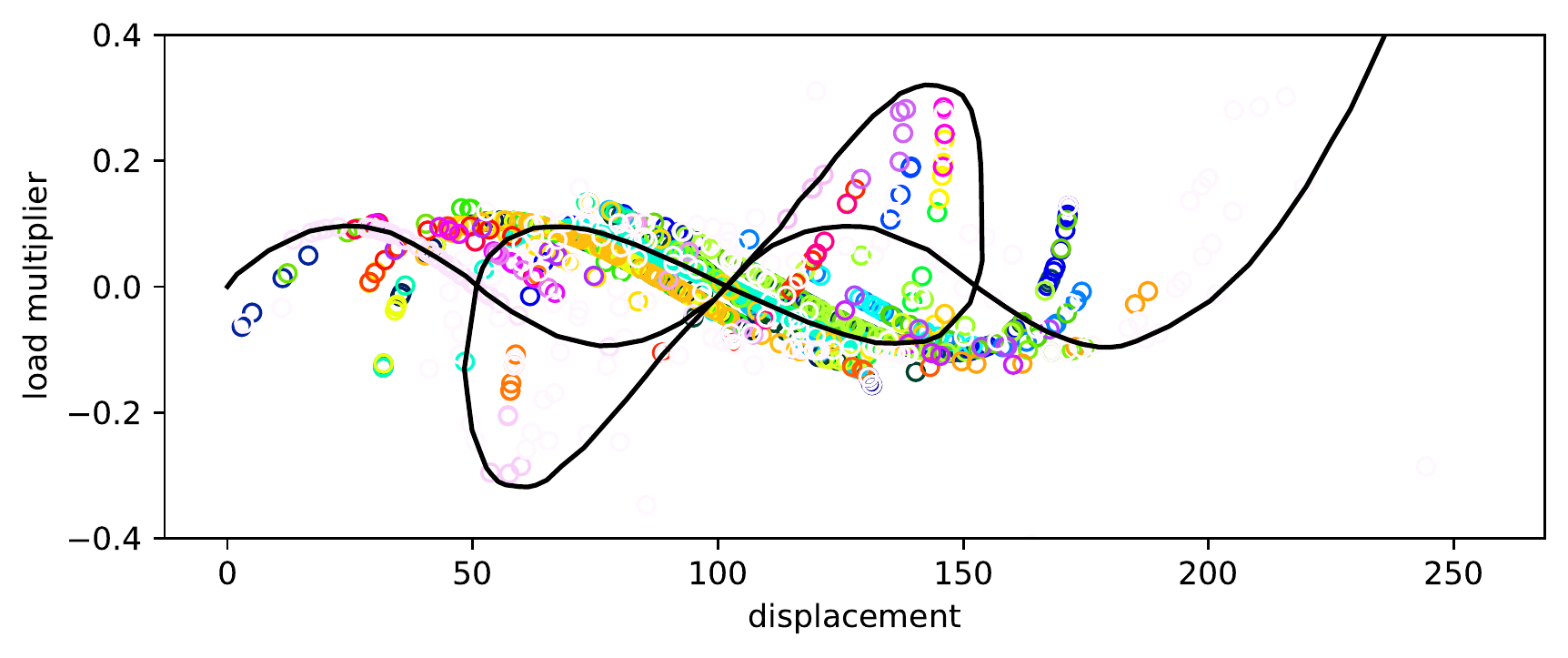}
        \caption{Clustering of solutions (obtained by DE/rand/1/bin algorithm) for identifying different feasible characteristics of equilibrium path (in different color circles).  \label{fig:3d_star_dbscan}}
\end{figure}

\textbf{Equilibrium path characteristics identification analysis.} A clustering analysis using the DBSCAN algorithm~\cite{ester1996density} was performed to analyze whether a few sets of solutions follow similar properties and whether they can be grouped together that identify particular equilibrium path characteristics. {\color{teal}The DBSCAN algorithm was chosen because its characteristics to scan through point cloud. The scanning of the neighborhood of the point cloud formed by solutions were also required for the equilibrium path characteristics identification, where the nearest connected points were expected to cluster together to represent particular path characteristic.}   Figure~\ref{fig:3d_star_dbscan} shows the clustering analysis results on DE/rand/1/bin solutions, and it does find a set of solutions that appear to follow particular equilibrium path characteristics. Hence, one can choose a cluster (a color) representing particular equilibrium path characteristics. 

One advantage of clustering analysis is that it offers several equilibrium paths compared to the standard civil engineering Arc-length method. However, solutions are not in order, and this approach produces too many equilibrium paths of varied unknown characteristics, especially many belonging to non-vertical displacement. Moreover, this analysis leaves with another issue, i.e., a cluster selection problem needed to be solved if one would like to choose only a single equilibrium path. Therefore, a more efficient method was required. In this research, a new algorithm called hypersphere search algorithms is presented for solving these class problems as accurately as possible.  

\subsubsection{Results of the \color{blue} adaptive search space decomposition method}
\label{sec:hypersphere_results}
Further analysis on the sixteen members truss structure problem was performed through the proposed hypersphere search algorithm. The algorithm started with a small domain of [\SI{-10}{\mm}, \SI{10}{\mm}] for the displacement variables and [\num{-0.2}, 0.2] for the load multiplier variable. First, a set of 5 trial solutions were {\color{red}evaluated} in step 1 pertaining to hypersphere initialization and first hypersphere center identification of this algorithm (cf. Section~\ref{sec:hypersphere_algo}). Then, a radius $r = \SI{5}{\mm}$ was used for the construction of the next hypersphere around the obtained center, and another set of 5 trial solutions were {\color{red}evaluated} to choose the next center for the next hypersphere, and this process continued until the control point reached the fixed maximum value of \SI{250}{\mm} or the maximum number of 1000 trials was reached.

\begin{figure}
        \centering
            \includegraphics[width=0.5\linewidth]{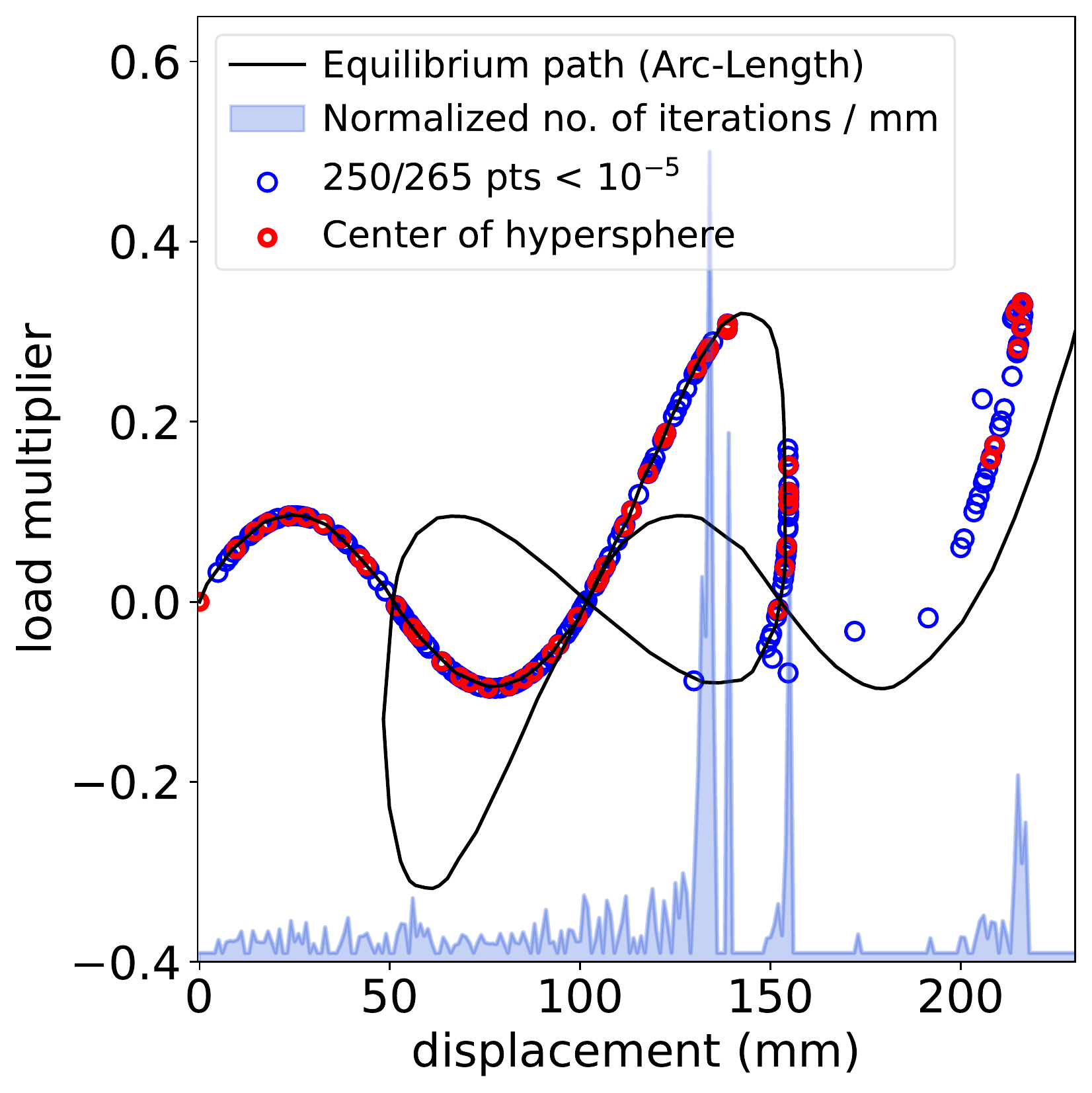}
        \caption{Quality of solutions produced by the hypersphere search algorithm with input radius $r = \SI{5}{\mm}$. Solutions in blue are trial solutions of respective hypersphere whose centers are indicated in red. {\color{blue} The computational effort in terms of normalized (between -0.4 and 0.5) total number of DE generations spent per mm displacement is shown in shaded blue area. The average computation effort is \num{8125} generations and max (i.e the peak) computation effort is \num{274168} generations. }
        \label{fig:3d_star_hypersphere_algo}}
\end{figure}

Figure~\ref{fig:3d_star_hypersphere_algo} shows the results of the hypersphere search algorithm. The obtained centers are shown in red, and all other optimal solutions are shown in cyan. This analysis has produced accurate solutions and a unique equilibrium path characterization for the vertical displacement. The only limitation of this method is the difficulty of this ill-posed nonlinear post-buckling analysis, where this algorithm finds it computationally challenging to search the next hypersphere where the structure has reached its breaking point (buckling point). That is, the radius $r = \SI{5}{\mm}$ for the hypersphere seems unable to find the next appropriate hypersphere due to sharp and sudden change in the search space (see a blue-shaded spike of the computational needs). This poses a highly challenging task for optimization algorithms.

\begin{figure}[h]
        \centering
            \includegraphics[width=0.45\linewidth]{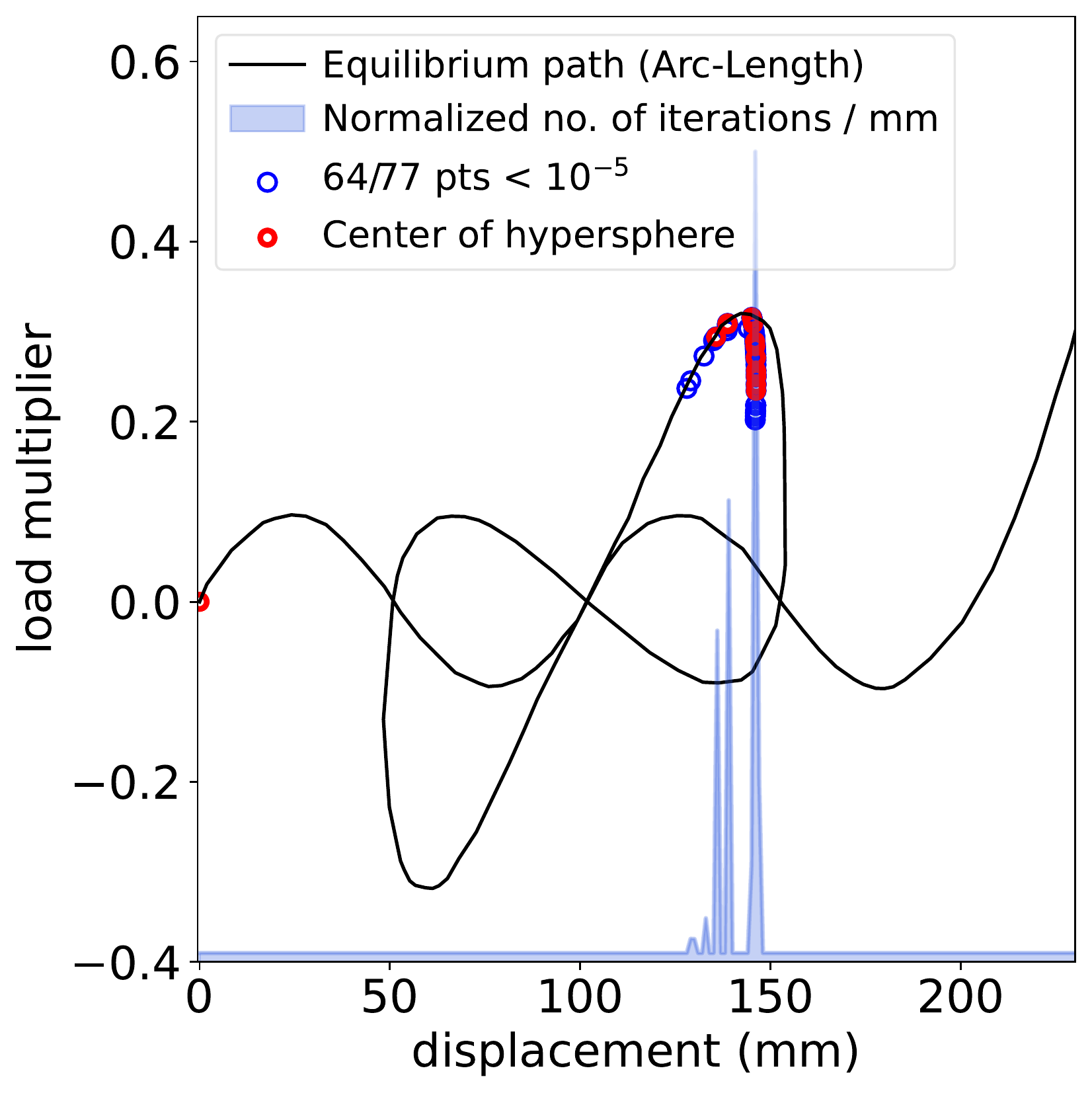}\hfill
            \includegraphics[width=0.45\linewidth]{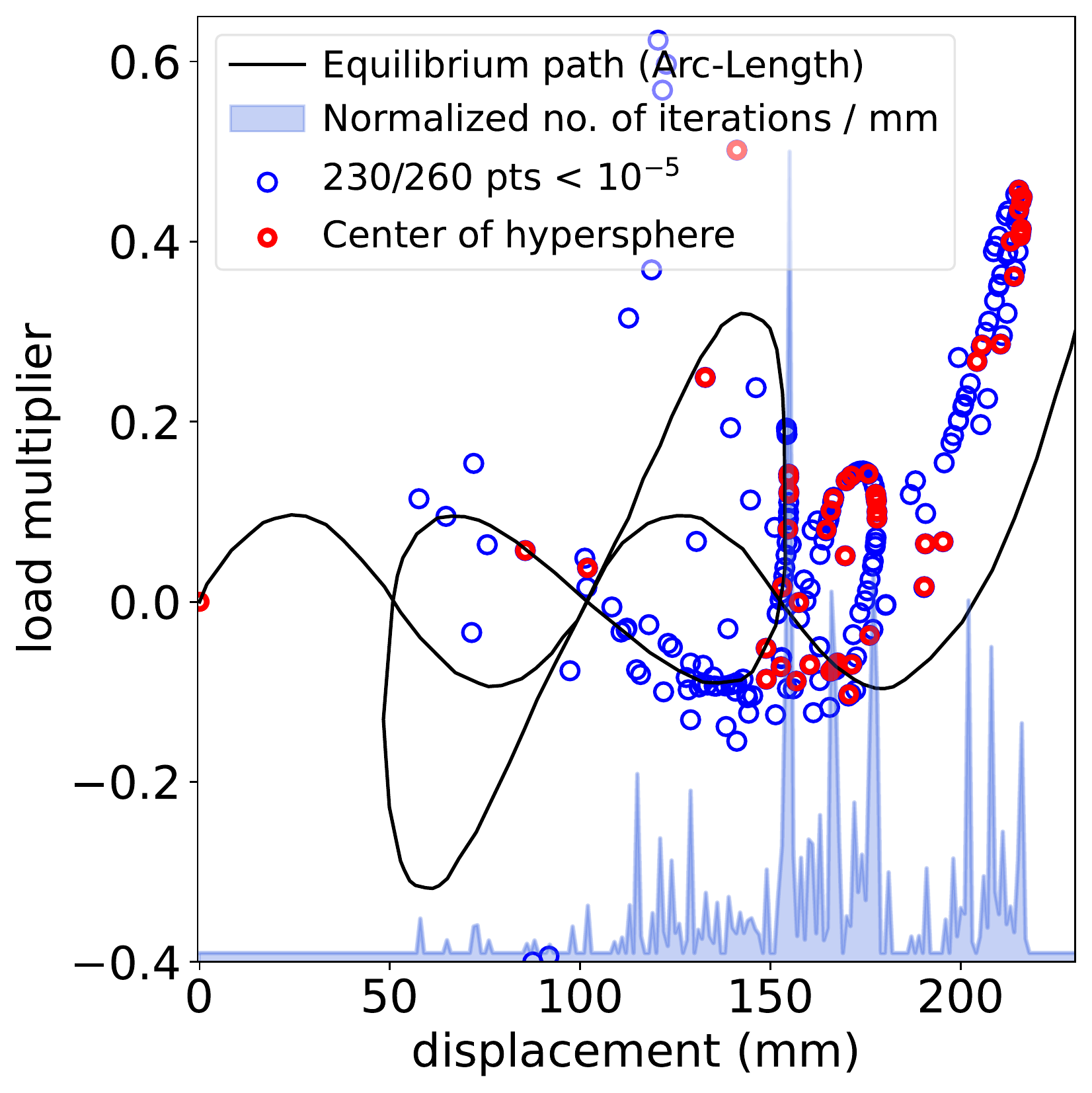}
        \caption{Sixteen-Member solutions produced by hypersphere search algorithm with input radius $r_{\text{new}} = r_{\text{prev}}\times 0.5$ (\textit{Left}) and with input radius $r_{\text{new}} = r_{\text{prev}} + \SI{5}{\mm}$ (\textit{Right}).{\color{blue} Reduction of the radius make algorithm stay at the buckling point and increase of the radius lead to non-vertical hypothesis equilibrium stages/path.} Solutions in blue are trial solutions of respective hypersphere whose centers are indicated in red. Both experiments in \textit{Left} and \textit{Right} were started from a ``center'' immediate behind the breaking points (see computational spike) obtained by the experiment in Figure~\ref{fig:3d_star_hypersphere_algo}.
        \label{fig:3d_star_hypersphere_algo_v2}}
\end{figure}

Another set of experiments was done to fine-tune user-defined hyperparameter radius $r$ of the hypersphere search algorithm. The first trial was to reduce to radius $r$ by a factor of half each time it was unable to find the next hypersphere (i.e.,  when it reached the breaking point). This was done by making radius $r$ adaptive to the iteration of hypersphere construction as per this expression: $r_{\text{new}} = r_{\text{prev}}\times 0.5$. The analysis in Figure~\ref{fig:3d_star_hypersphere_algo_v2}(left) shows that the algorithm can follow the sharp breaking path; however, it is computationally challenging and slow. {\color{blue} This is because eventually, the hyperspheres become small and infeasible for the algorithm to continue. The computational effort of this stage is shown by the spikes (peak values) of computational effort in Figures~\ref{fig:3d_star_hypersphere_algo} and \ref{fig:3d_star_hypersphere_algo_v2}.}

Another trial was an additive approach where radius $r$ was gradually increased by a value of 5mm, i.e., $r_{\text{new}} = r_{\text{prev}} + 5$. Figure~\ref{fig:3d_star_hypersphere_algo_v2}(right) shows the algorithm was able to escape the breaking point, but it deviated from finding an exact unique path, and it jumped to other feasible solutions {\color{teal}(belonging to non-vertical load)} like the many alternatives found in clustering analysis in Section~\ref{sec:star_search_domain_anal}.  

\subsection{\color{blue}Benchmark 3: Twenty four member shallow truss structure}
\label{sec:star_new}
{\color{blue}
The third benchmark, represented in Figure~\ref{fig:3d_star_new}(left side), is a space truss structure composed of 24 trusses, whose axial rigidity is \SI{960.5}{\kN}, seven free nodes, and 21 free displacements. The permanent load comprises a vertical force of 50N applied on the central node and six vertical forces applied on the other free nodes of the structure. The structure has been analysed in~\cite{rezaiee2011automatic} using a dynamic relaxation method. Therefore, the load path curve of the structure, referring to the vertical displacement of the central node, is known in the range from \SI{0}{\mm} to \SI{50}{\mm}. Figure~\ref{fig:3d_star_new} (right side) shows the load-deflection curve of the central node and the solutions provided by the adopted optimization procedure. The two critical points are located approximately at \SI{8}{\mm} and \SI{30}{\mm}. Between these displacement values, the structural response is characterized by an evident snap-through mechanism. 

The hypershpere method for the adaptive search space decomposition together with DE algorithms provided a good-densed distribution of solutions in most of the investigated range of displacements, including the two critical zones where a limited number of iterations were required for the convergence. The proposed methodology was also capable of describing the snap-trough path with a significant number of solutions. However, a significantly higher number of iterations were required in this response stage. In addition, from the comparison, it is possible to observe that the solutions provided by the proposed methodology follow the baseline with a good level of accuracy within the entire investigated range of displacements. Finally, Figure~\ref{fig:3d_star_new_shapes} shows different deformed shapes of the structure at different displacement magnitudes.     


\begin{figure}[H]
        \centering
             \includegraphics[width=0.45\linewidth]{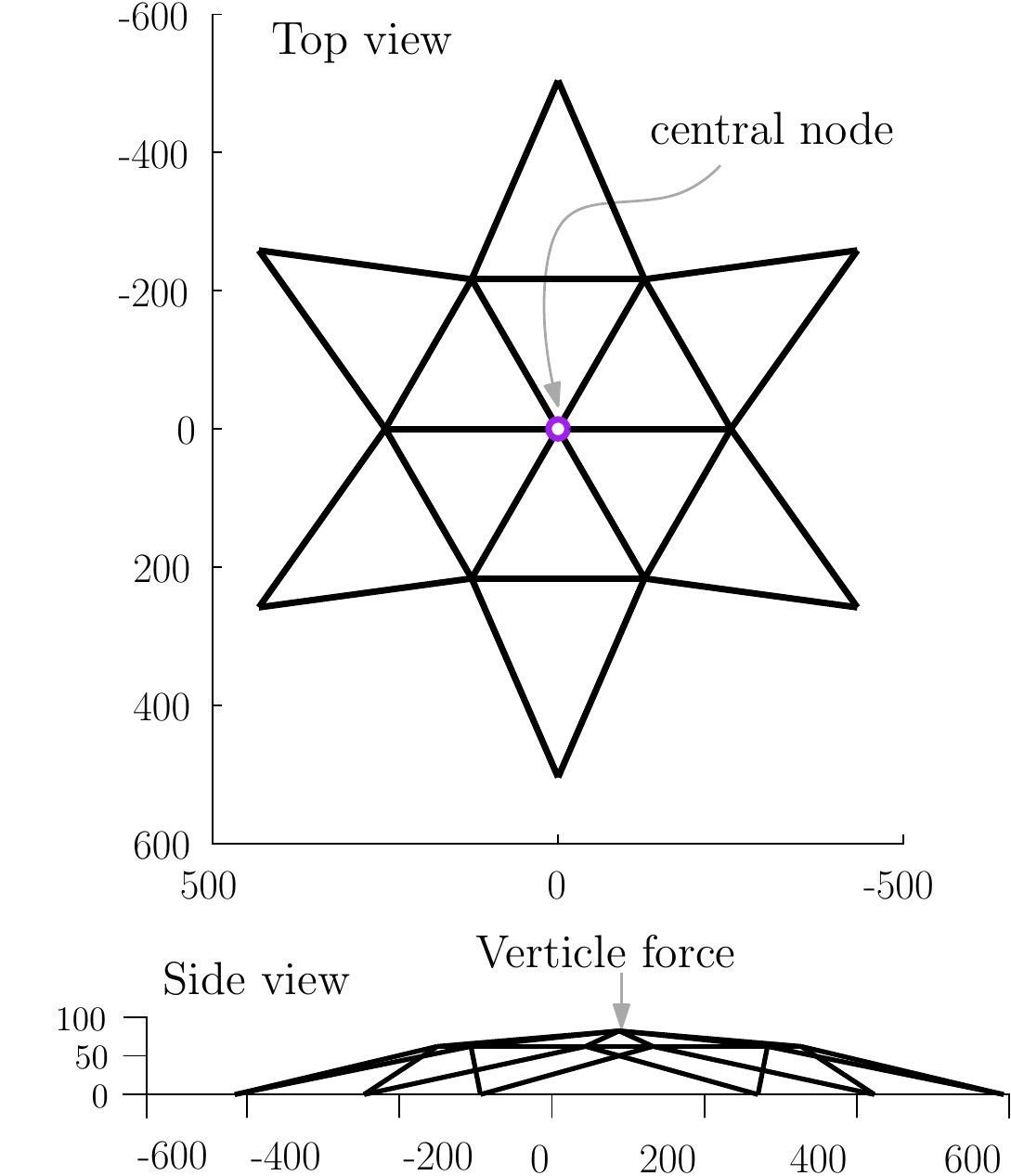} \quad
            \includegraphics[width=0.5\linewidth]{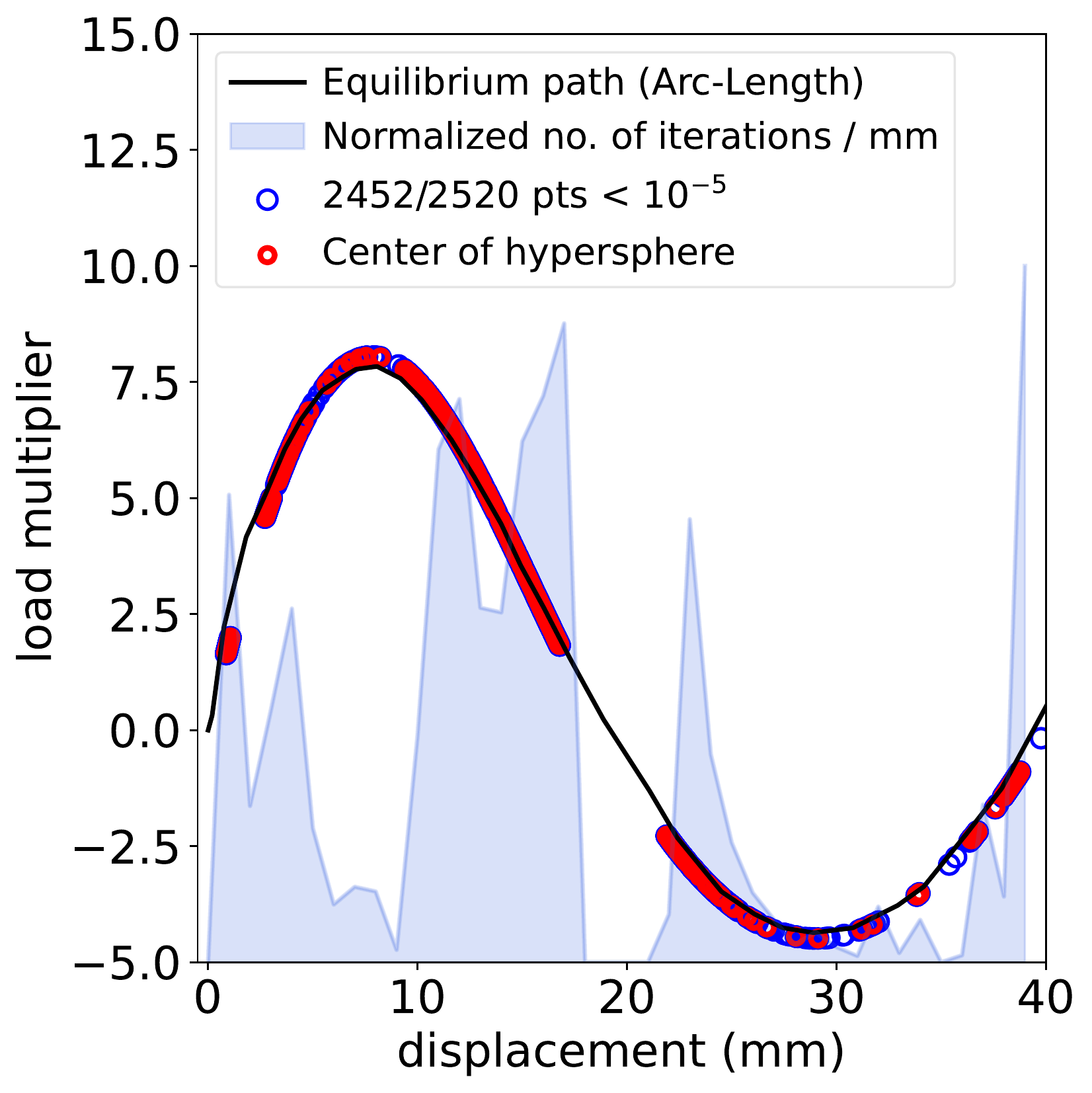}
        \caption{\color{blue} Twenty Four Member Shallow Truss Structure. \textit{Left:} Top View and Side View of the structure. The central node of the structure is indicated by a lighter color dot in Top view and a light color arrow pointer in Side view. The arrow on Side view indicates a vertical downward external force of \SI{960.5}{\kN} applied on the central node. \textit{Right:} The results of adaptive search space decomposition method applied on this benchmark. Blue circles are optimal solution, red circles are centers of the hyperspheres of the adaptive search space decomposition method. Shaded curve area is the computational effort in terms of normalized (between -5 and 10) sum of the number of generations spent of searching solutions at a hypersphere center per mm displacement. The average computation effort is \num{24468} generations and max (i.e the pick) computation effort is \num{37369}. 
        \label{fig:3d_star_new}}
\end{figure}

\begin{figure}[H]
        \centering
        \includegraphics[width=0.45\linewidth]{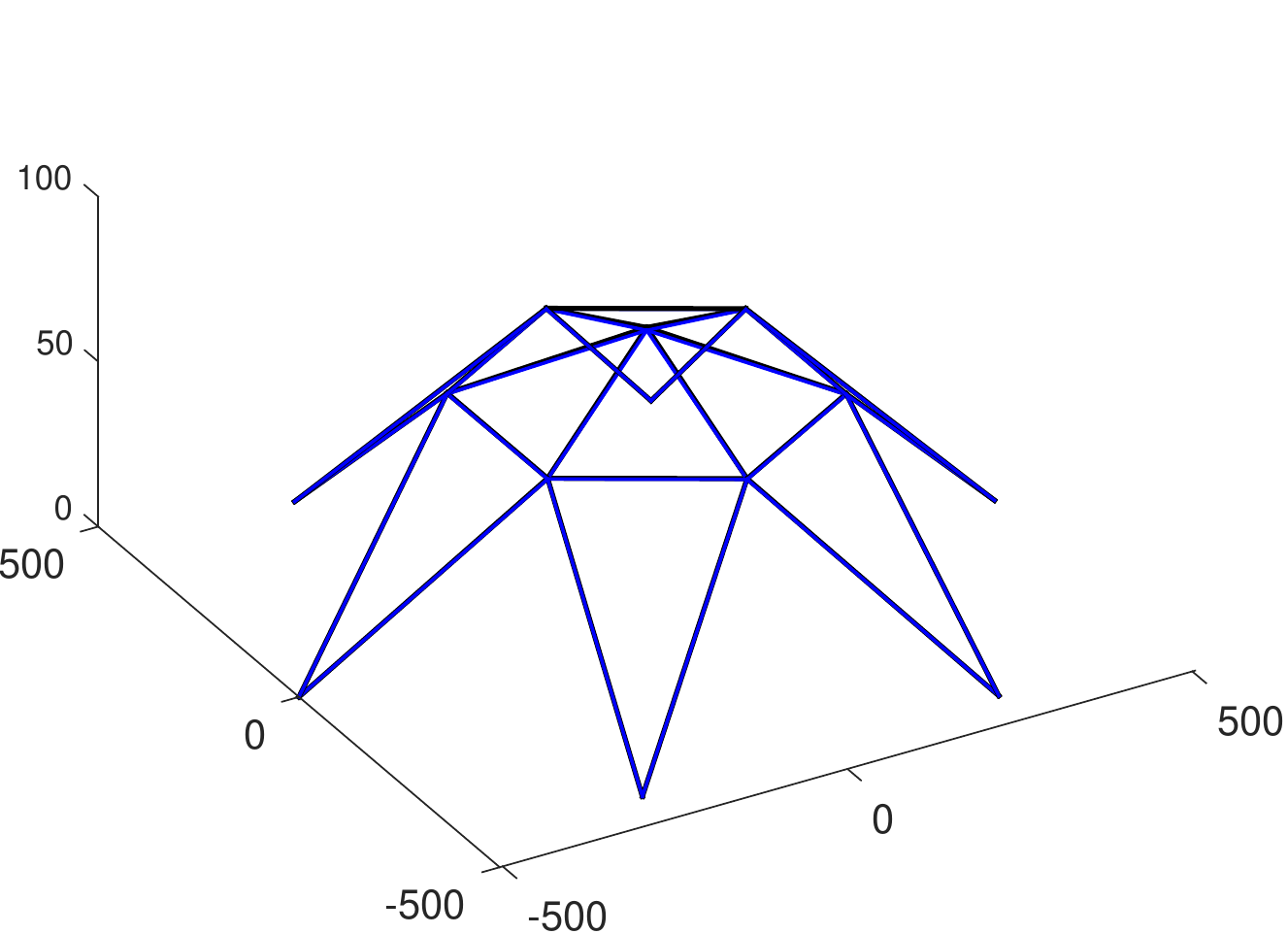} 
        \includegraphics[width=0.45\linewidth]{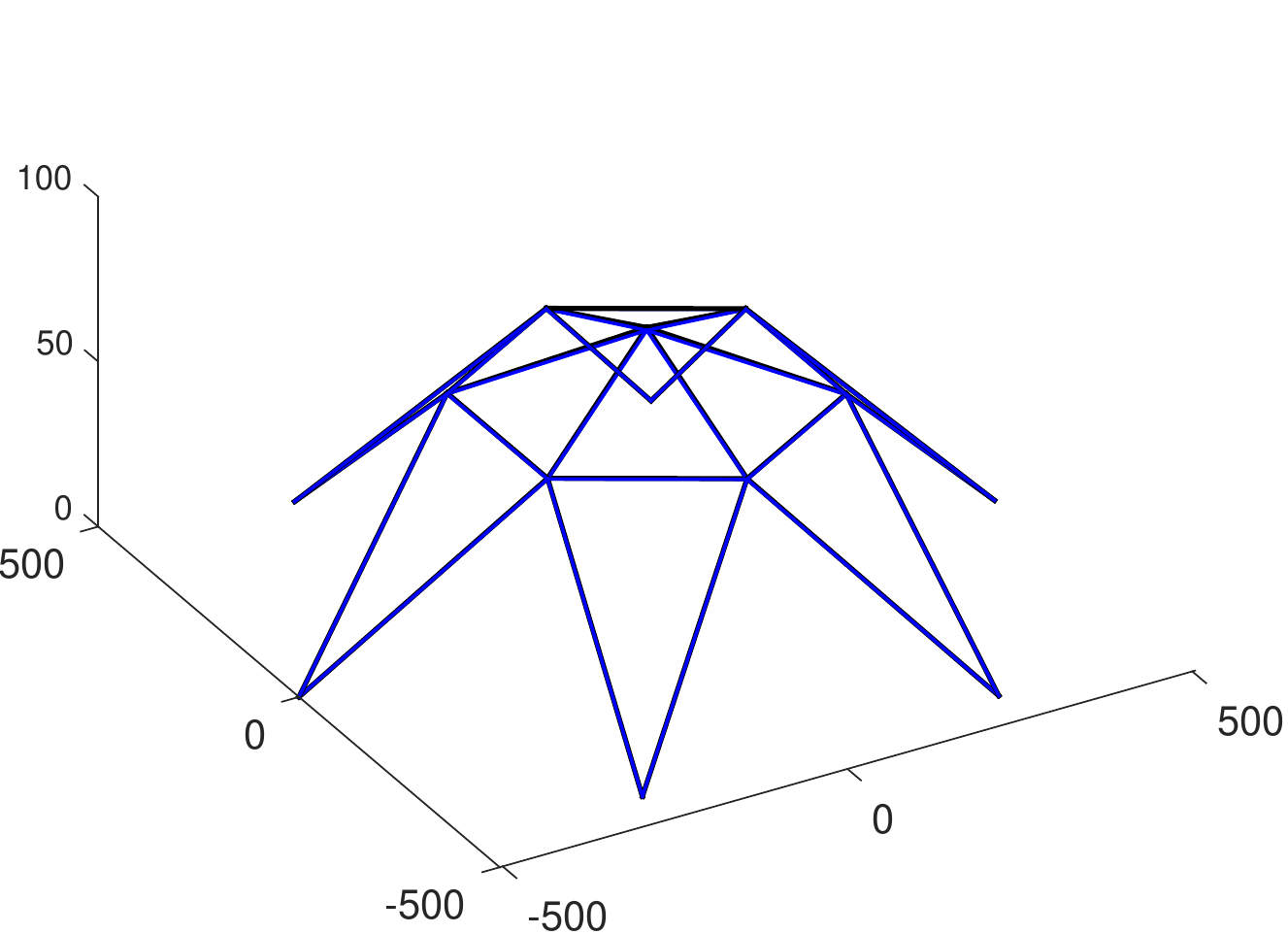} 
        displacement 0.8 mm and load 1.64 kN \quad\quad  displacement 3.3 mm and load 4.99 kN
        \includegraphics[width=0.45\linewidth]{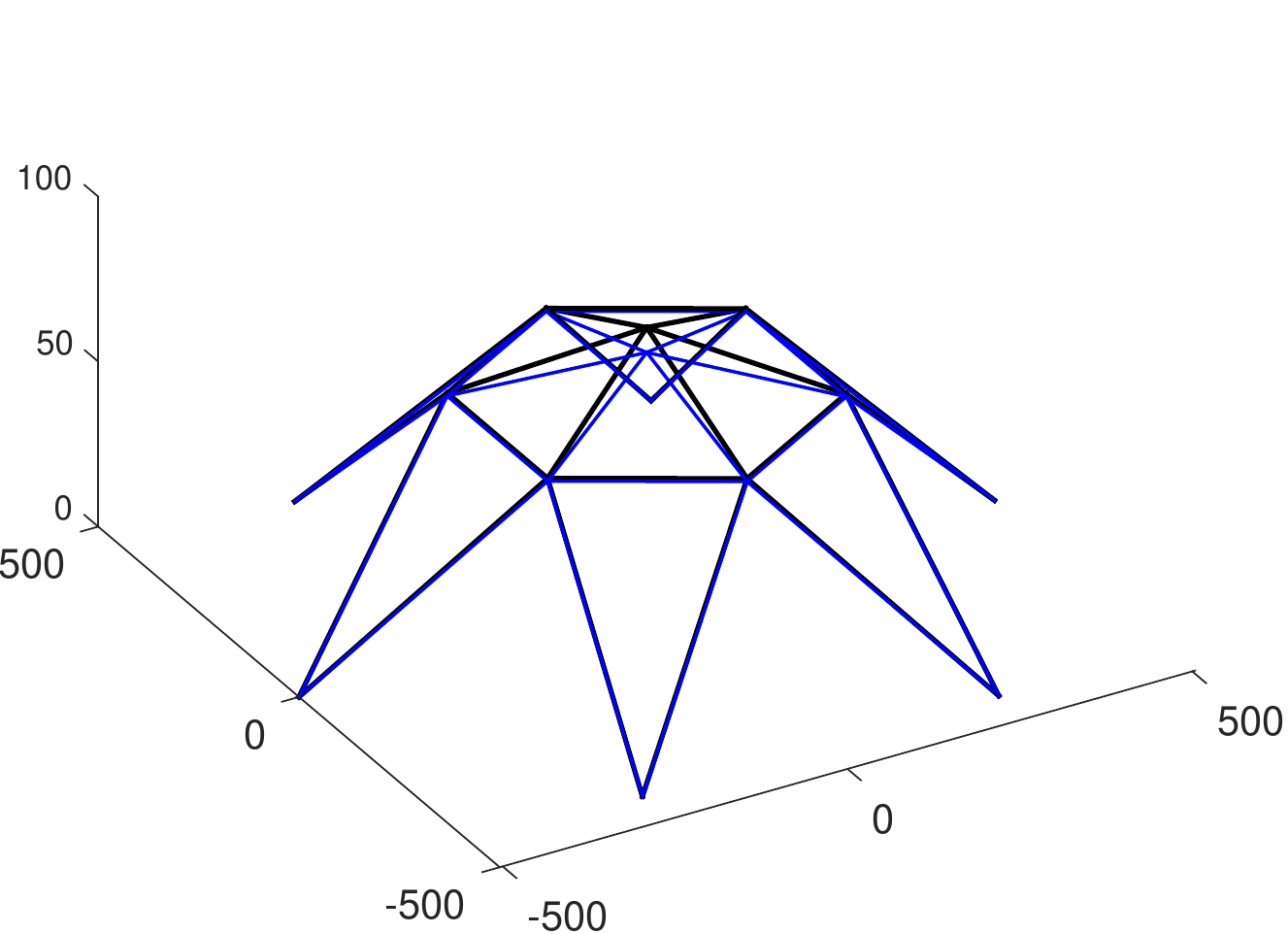} 
        \includegraphics[width=0.45\linewidth]{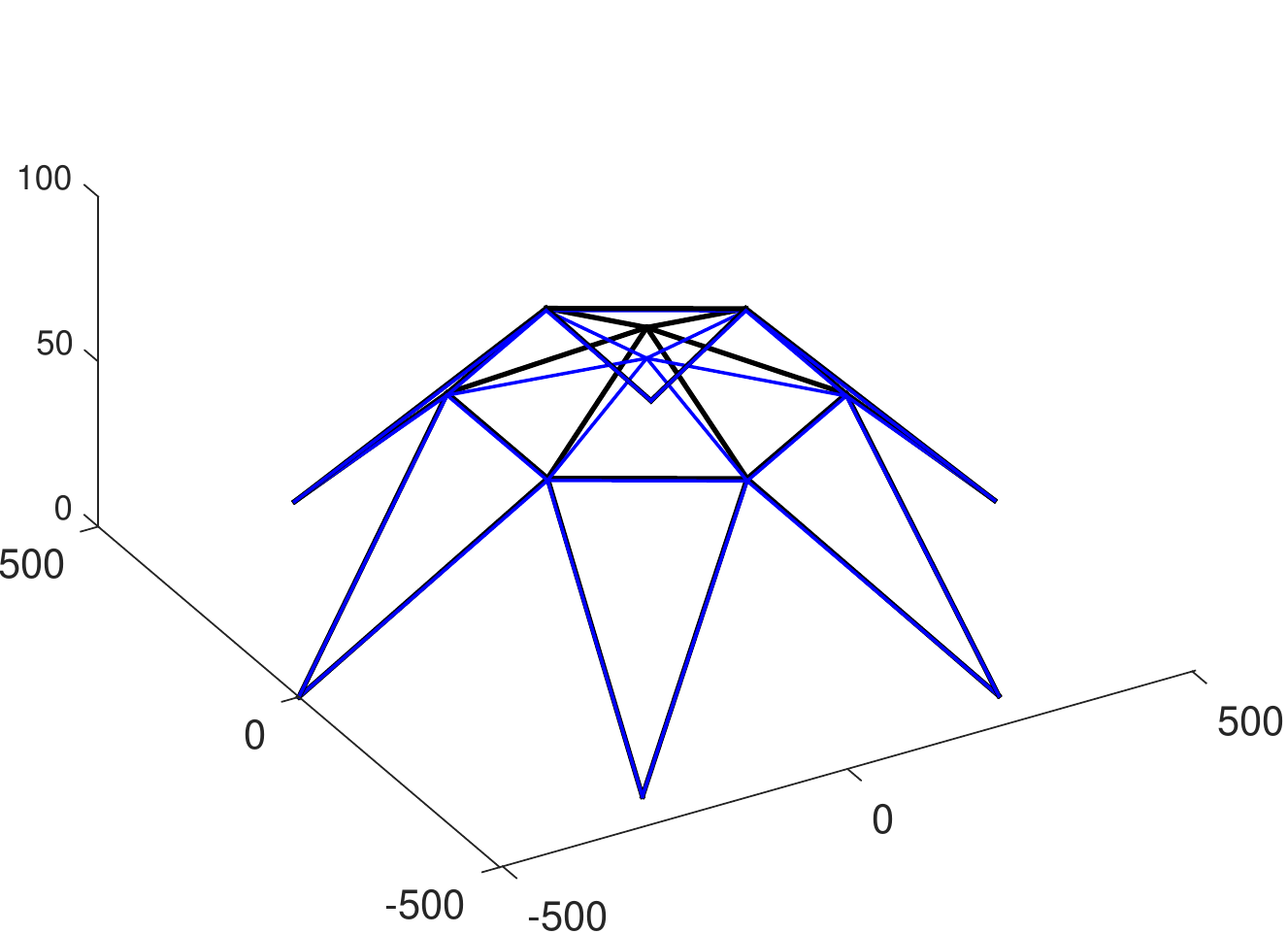} 
        displacement 7.5 mm and load 8.04 kN \quad\quad  displacement 9.3 mm and load 7.78 kN
        \includegraphics[width=0.45\linewidth]{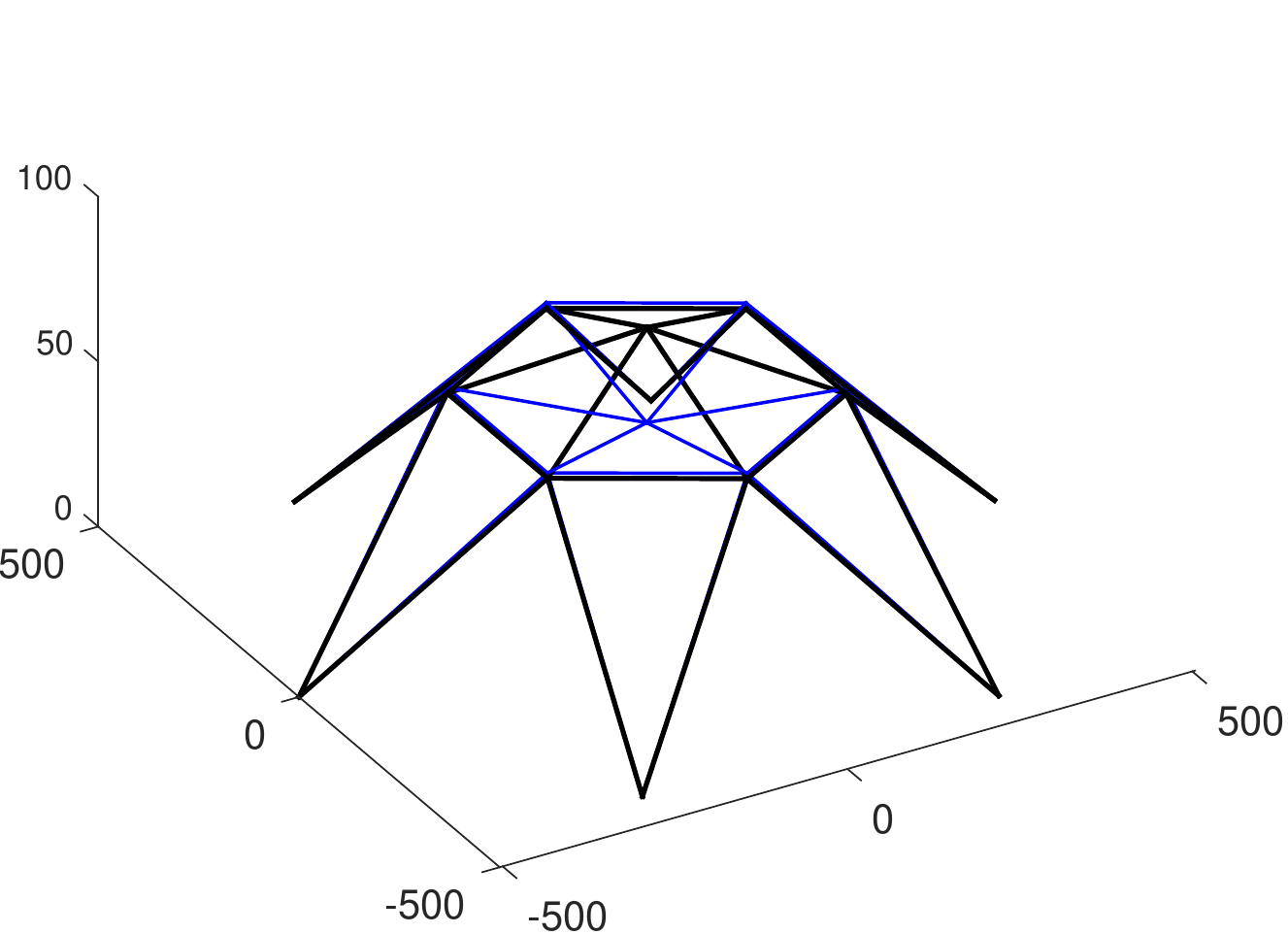} 
        \includegraphics[width=0.45\linewidth]{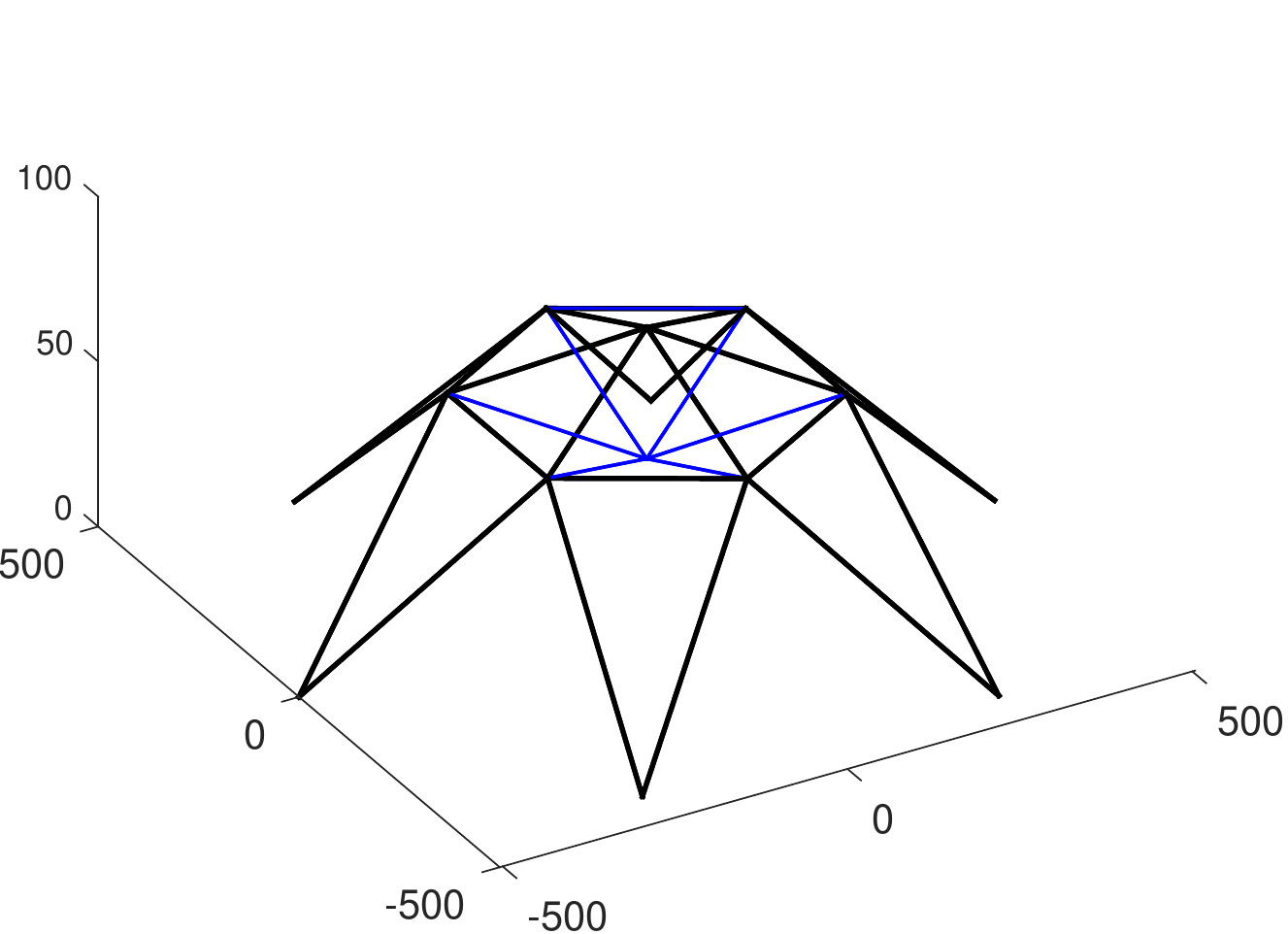} 
        displacement 28.5 mm and load -7.49 kN \quad\quad  displacement 40 mm and load -0.17 kN
        \caption{\color{blue} Deformed shapes and equilibrium stages of twenty four member shallow truss structure member space truss structure.
        \label{fig:3d_star_new_shapes}}
\end{figure}
}

\subsection{Test Problem: \color{blue}spatial reticular beam structure}
\label{sec:bridge}
\begin{figure}
	\centering
	Side view (2D view)\hfill Top view (3D view)\\[2pt]
	
	\includegraphics[width=0.45\textwidth]{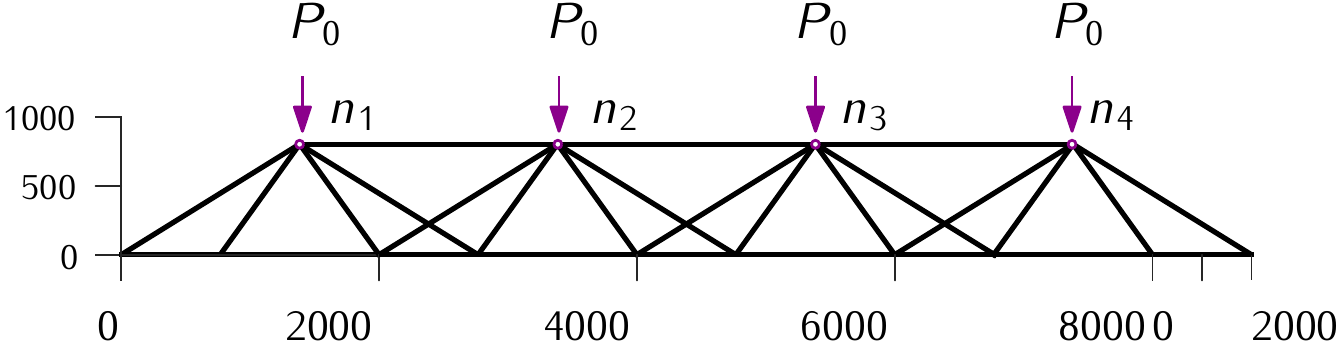}~\quad
    \includegraphics[width=0.45\textwidth]{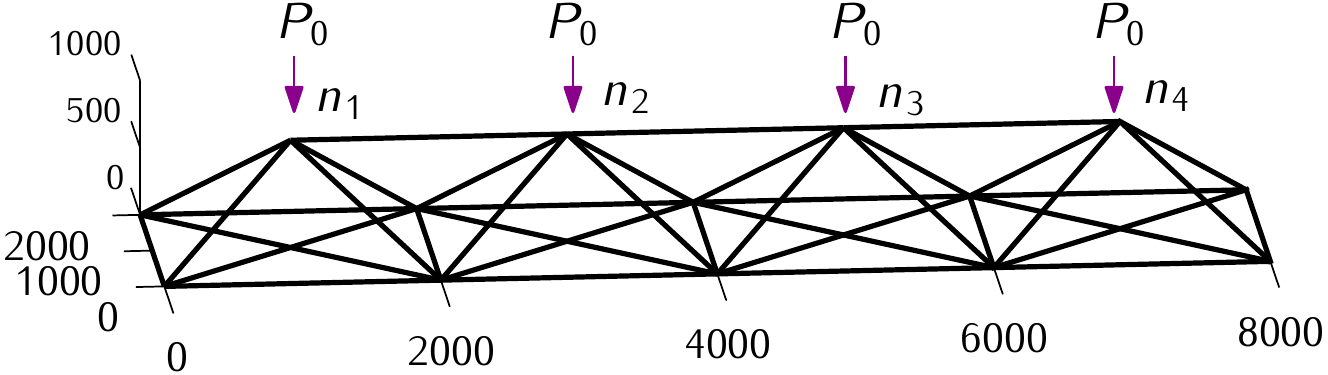}
    Undeformed shape
    
	\includegraphics[width=0.45\textwidth]{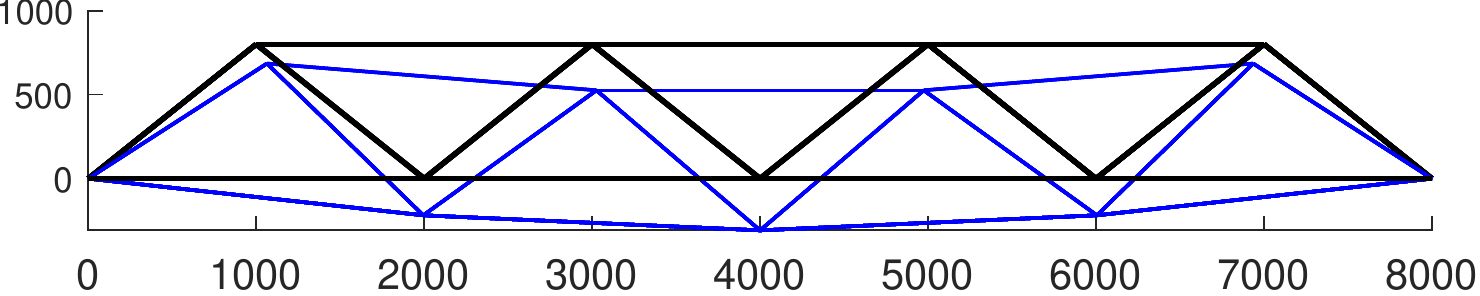}~\quad
	\includegraphics[width=0.45\textwidth]{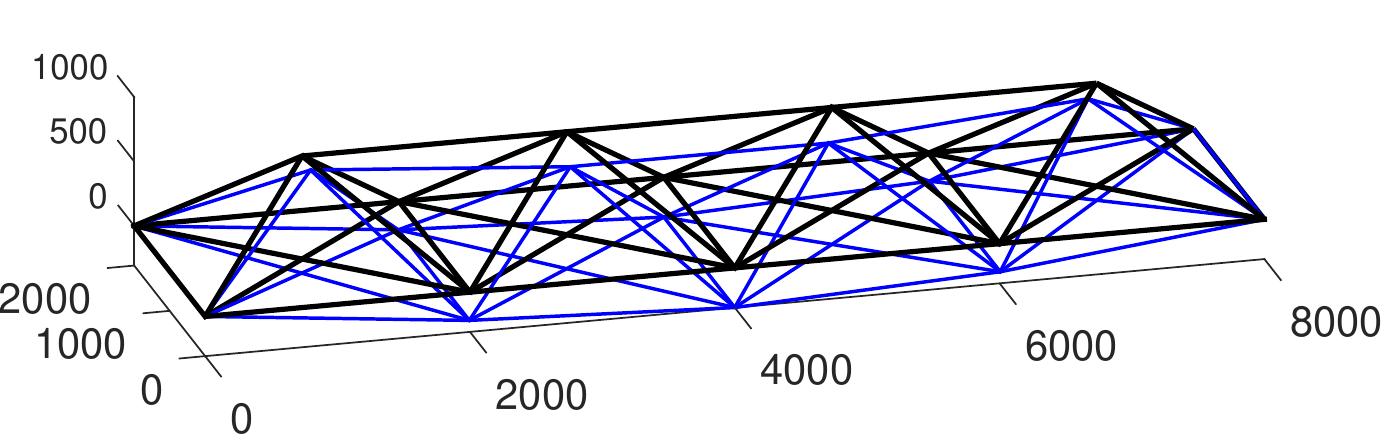}
    Deformed shape A
    
	\includegraphics[width=0.45\textwidth]{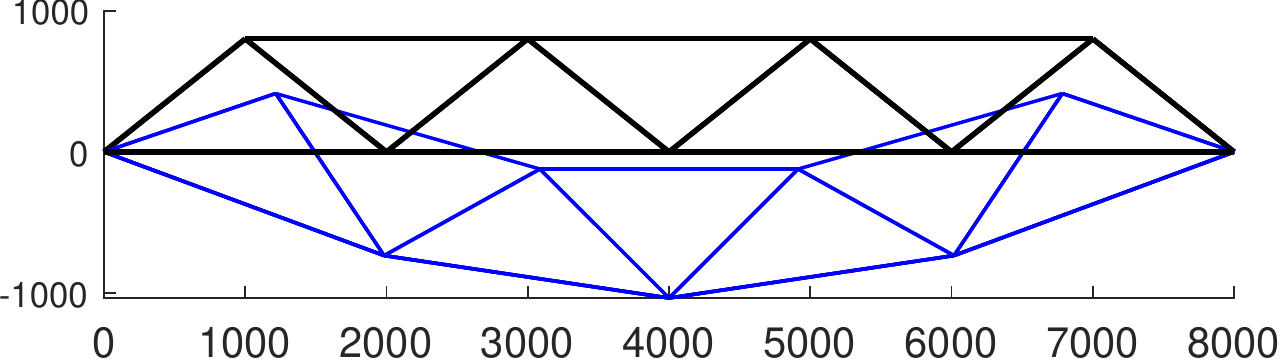}~\quad
	\includegraphics[width=0.45\textwidth]{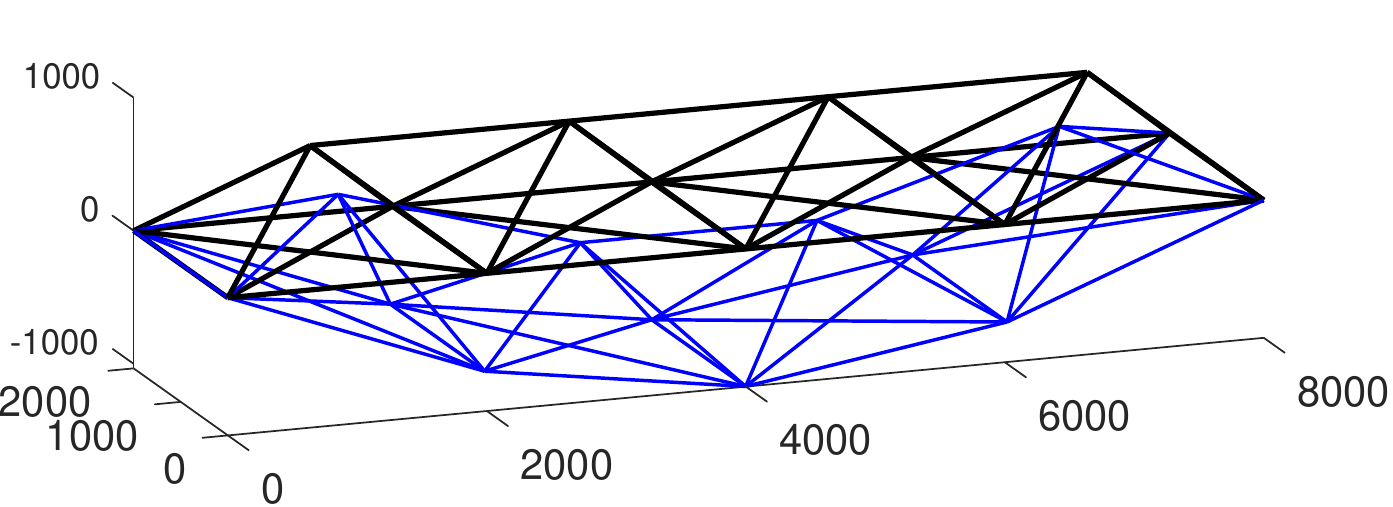}
    Deformed shape B
    
	\includegraphics[width=0.45\textwidth]{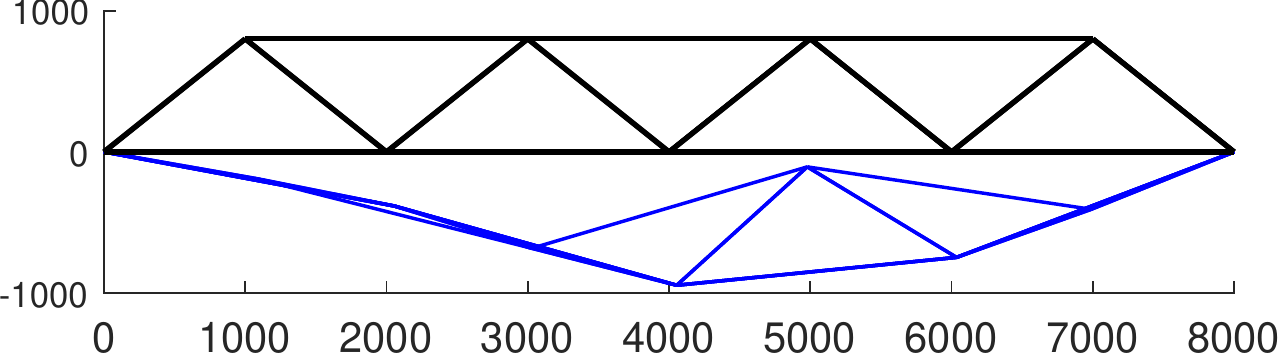}~\quad
	\includegraphics[width=0.45\textwidth]{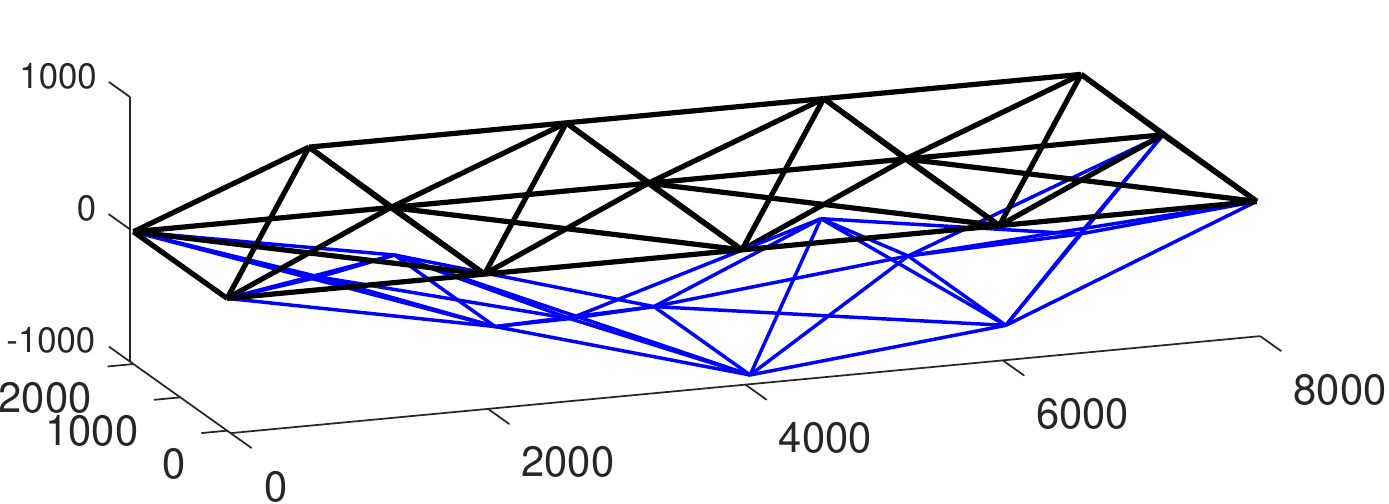}
    Deformed shape C
    
	\includegraphics[width=0.45\textwidth]{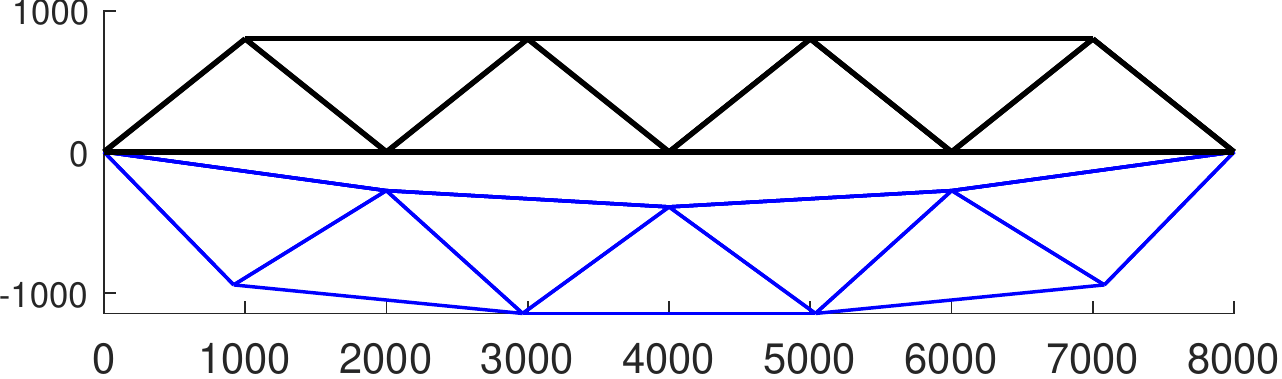}~\quad
	\includegraphics[width=0.45\textwidth]{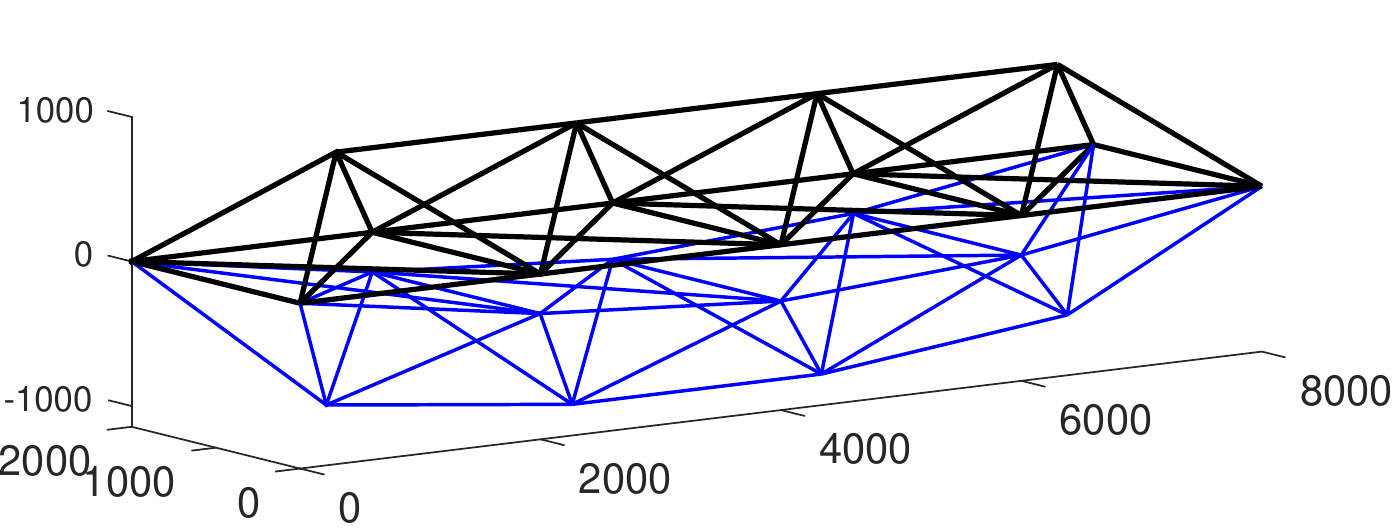}
	Deformed shape D
	
	\caption{Equilibrium stages of 3D reticular beam structure optimized using DE algorithm version DE/rand/1/bin. Lighter color arrows (on Side view and Top View) indicates the application of the vertical downward force ($P_0 = \SI{100}{\kN}$) applied on the nodes $n_1, n_2, n_3,$ and $n_4$ of 3D reticular beam structure. The undeformed (original) structure shapes are drawn with black lines and deformed structure shapes are drawn with blue lines. The structure shapes 2D view is in the left column and their 3D view is in the right column. All deformed shapes presented achieved an accuracy of \num{1e-5} on the objective function Equation~\eqref{eq:main_objective}. Row 1 is original structure and rows 2 (Shape A), 3 (Shape B), 4 (Shape C), 5 (Shape D) are deformed shapes for respective applied external lead multipliers values 25.6302, 158.6989, 13.5174, 55.8157. The deformed shapes A, B, C, and D respectively has their control point displacements 274.2566, 921.2034, 1482.2, and 1943.6.
	\label{fig:bridge_view}}
\end{figure}

This Section presents the results of the analysis on a medium space truss structure representative of a 3D steel beam supporting a roof system. The structure has a span of \SI{8000}{\mm}, a width of \SI{2000}{\mm}, and a height of \SI{750}{\mm}, arranged with the geometrical layout shown in Figure~\ref{fig:bridge_view}(row 1). The beam is composed of 33 trusses, 10 free nodes, and 30 DoF. Each truss has \SI{2500}{\mm^2} of cross-section area and \SI{200000}{\MPa} of Young’s modulus. Four simple supports are considered at the two ends of the beam, and the external load is represented by four vertical forces, equal to \SI{100}{\kN}, applied at the beam top nodes.

On this problem, DE/rand/1/bin and DE/best/2/bin algorithms were applied as they were the most suitable algorithms for this class of structure analysis problem in the benchmark problem optimization. The goal was to find a few feasible solutions to approximate an equilibrium path containing pre-and post-buckling stages of a 3D steel beam. Hence, the variable domain setup was [\SI{0}{\mm}, \SI{2500}{\mm}] for the displacement and [1, 300] for the load multiplier. In this setting, both DE algorithms versions were applied with stopping criteria of 1 million iterations or convergence accuracy of \num{1e-5}. 
Several instances of both versions of DE reached an accuracy of \num{1e-5}. However, when filtering solutions, DE/rand/1/bin was found to be producing a greater number of solutions that were accurate compared to DE/best/2/bin. Hence, solutions obtained by DE/rand/1/bin are shown in Figure~\ref{fig:bridge_view}.  In Figure~\ref{fig:bridge_view}, rows 2 to 5 show deformed shapes of four equilibrated and feasible solutions. 

\begin{figure}
	\centering
	\includegraphics[width=0.8\textwidth]{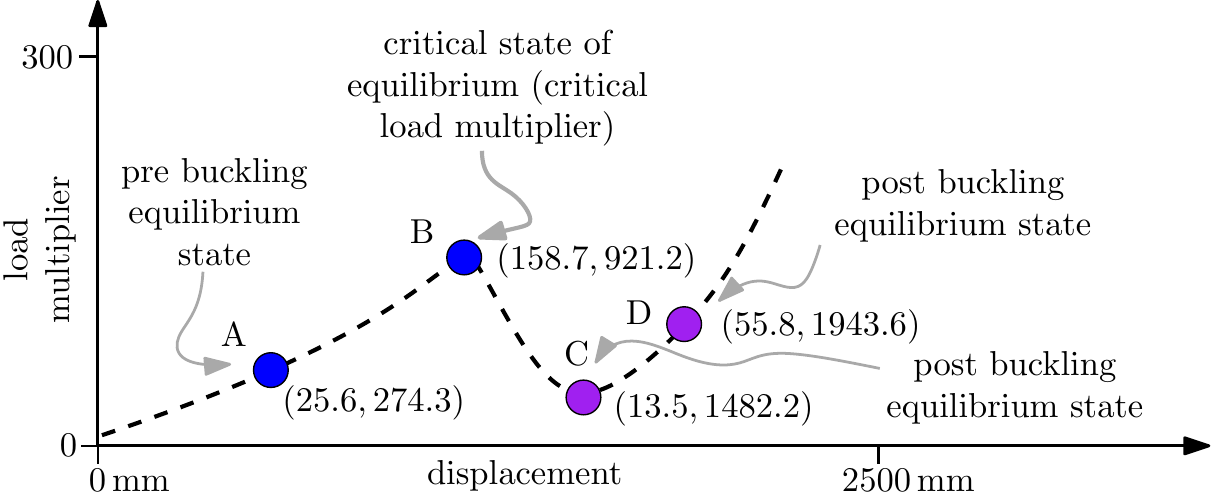}
	\caption{Pre- and post buckling equilibrium states of the 3D reticular beam. These states are obtained for a multistep analysis mentioned in Figure~\ref{fig:Multi_Step}. The dotted line is the hypothesis equilibrium path connecting the obtained states for this test problem. \label{fig:bridge_points}}
\end{figure}

A rough representation of the equilibrium path is represented in Figure~\ref{fig:bridge_points} in terms of maximum vertical displacement of the structure and the load multiplier to characterize these solutions. The dashed lines connecting the solutions represent qualitatively the equilibrium path between these solutions which can be drawn by considering a refined set of solutions. However, despite the limited number of solutions considered, it was possible to highlight the main characteristics of the pre- and post-buckling behavior of the structure. 

More in detail, it is possible to observe that solutions A and B belong to the stable equilibrium (pre-buckling) region, where the increase of deformations corresponds to the increase of the external loads and the deformed shape is symmetric (cf. Figure~\ref{fig:bridge_view} and \ref{fig:bridge_points}). On the other hand, solution C shows a buckling behavior (cf.  Figure~\ref{fig:bridge_view} and \ref{fig:bridge_points}) involving the truss elements of the right part of the beam, leading to the decrease of the load multiplier and the loss of the symmetry of the deformed shape. Finally, Solution B is characterized by the global buckling with the overturning of the structure (cf. Figure~\ref{fig:bridge_view} and \ref{fig:bridge_points}). After this point, the system behaves as a catenary rather than a beam, and the equilibrium is again stable.          

\section{Discussions}
\label{sec:discussion}
This paper presents a novel methodology based on gradient-free optimization algorithms for assessing the {\color{red} nonlinear} structural response of space truss structures subjected to large displacements, {\color{red} accounting for geometrical nonlinearities.} {\color{blue} This methodology allows for evaluating the load path of the structure, including stable and unstable equilibrium stages, thus characterizing its post-buckling behavior and assessing the vulnerability to progressive collapses.} Compared to the standard Newton-Raphson procedures, the proposed methodology does not require the assemblage and update of the global stiffness and geometric {\color{red}matrices} (Section~\ref{sec:opt_problem}). Moreover, it is not exposed to the potential numerical issues affecting the Newton-Raphson and arc-length methods in correspondence to the critical points, where the equilibrium changes from the stable to the unstable.

This research applied four classes of optimization algorithms: deterministic optimization algorithm (DIRECT), single solution-based algorithm (Simulated Annealing), swarm inspired algorithms (ABC, ACO, and PSO), and an evolutionary algorithm (DE) (Section~\ref{sec:opt_algo}). {\color{teal} The DE algorithm was found to be better at converging among all algorithms from our initial trial of three benchmark problems. Hence, the scope of this work mainly relied on the DE algorithm as it was competitively the most successful in solving this class of problems.} {\color{blue} Finally, this study developed a multistep analysis algorithm, the hypersphere algorithm for the adaptive decomposition search space to support gradient free algorithms effectively describe the structure load path.}    

{\color{teal} Three benchmarks of 3D truss structures and one test case representative of a medium-sized structural problem were tackled} as optimization problems. {\color{blue} The optimization problem associated with the global structural equilibrium (Section~\ref{sec:opt_problem}) can be classified as a nonlinear, multimodel, unconstrained minimization optimization problem, which presents significant challenges to optimization algorithms.} {\color{teal}It was highly challenging to solve} as variables were highly sensitive to each other due to their physical interactions (see Section~\ref{sec:results}). Moreover, this optimization problem has a significantly large number of local minima. Additionally, this problem has multiple global minima, but the global minima were strongly linked to the load multiplier variable and the direction of loadings (downward forces on a structure), making this problem a nonlinear, multimodal, unconstrained optimization problem. {\color{blue} Based on the results of numerical experiments conducted in Section 3, it was possible to draw the following  conclusions:}  

\begin{itemize}
    \item The DIRECT algorithm was able to find solutions for the first benchmark but could not explore any other solutions except when the search space was manually partitioned to help it find other solutions and when the search space was highly restricted after domain analyses (cf. Figures~\ref{fig:eight_algo_point} and~\ref{fig:eight_disc_point}). Apart from that, DIRECT struggled to find any other solutions to benchmark problems. This is because the hyperspace became ill-defined or search space too complex for DIRECT to work properly. A similar observation is applied to the Simulated annealing algorithm, which was able to solve the first {\color{teal}problem, a small benchmark problem consisting of four variables,} but SA could not effectively solve the 16 variables sixteen-member problem.

    \item ACO demonstrated a tendency to find a solution only at the center of the domain. A close examination suggests that the ACO uses Gaussian distribution to generate new solutions in its iterations, and that may lead the solutions following the denser region of the search space (cf. Figures~\ref{fig:eight_algo_point} and ~\ref{fig:3d_star_algo_points}). However, this observation is merely a hypothesis given the nature of this algorithm and this problem as observed during the experimentation. 
    
    \item Swarm-inspired PSO algorithms were able to solve the first benchmark, but they had relatively poor performance on the second benchmark (cf. Figures~\ref{fig:eight_algo_point} and ~\ref{fig:3d_star_algo_points}). The DE versions were the best performing algorithms for this class of ill-posed problems. 
    
    \item Despite being the best performing, DE versions were unable to find the expected number of accurate solutions on the equilibrium path (Figure~\ref{fig:3d_star_full_disc_domain_anal}). Only after a domain decomposition analysis to better inform the algorithms about the domain range of variables, DE versions were able to find a high percentage of accurate solutions  (Figure~\ref{fig:3d_star_full_disc_domain_anal} and  Figure~\ref{fig:3d_star_algo_points}).

    \item Among DE/best/2/bin and DE/rand/1/bin, the DE/rand/1/bin was found to be best performing for such challenging optimization problems; the convergence profiles of both show that DE/best/2/bin tended to fall into local minim more often than the DE/rand/1/bin. This was attributed to solutions following the local best solution leading to local minima (Figure~\ref{fig:3d_star_convergence}). 

    \item On the test problem, the 3D reticular beam, which is a {\color{teal} real medium-size structure}, DE/rand/1/bin, like in the case of the second benchmark, outperformed DE/best/2/bin. It was able to find more optimal solutions for the test problem 3D reticular beam. The results of DE/rand/1/bin presented in Section~\ref{sec:bridge} have high precision and produce solutions for both pre- and post buckling stages (Figure~\ref{fig:bridge_view}).  These solutions allowed producing {\color{teal} hypothesizing the} equilibrium path (Figure~\ref{fig:bridge_points}).
    
    \item The proposed hypersphere algorithm, i.e. adaptive search space decomposition method, incrementally constructed hyperspheres using the knowledge from the previous hyperspheres (Section~\ref{sec:hypersphere_algo}). This {\color{teal} method helped DE find} solutions accurately, and its performance was comparable to the arc-length method used in civil engineering for the equilibrium path analysis (Figures~\ref{fig:3d_star_hypersphere_algo} and \ref{fig:3d_star_new}). 
    
    \item The main challenges with the proposed adaptive search space decomposition method were to tune the value of the user-defined hyperparameter radius of the hyperspheres, which included whether to fix it a single value for all iterations or make it adaptive to iterations (Section~\ref{sec:hypersphere_results}). Additionally, this algorithm suffers from an obvious issue concerning sudden and sharp change in the hypersphere's domain at the breaking/buckling point of the structure. Since the buckling point changed the domain sharply and suddenly, this algorithm struggled to find the next hypersphere (Figure~\ref{fig:3d_star_hypersphere_algo_v2}). 
\end{itemize}

{\color{blue} In summary, the main outcome of this research is to demonstrate that the use of heuristic-based optimization algorithms can be adopted for solving nonlinear structural problems, in particular for assessing the geometrically-nonlinear response of space truss structures, which are particularly prone to show complex, unstable equilibrium stages and snap-through mechanisms. This outcome was achieved by successfully solving three} benchmark problems to a high degree of accuracy and provided several optimal solutions for the test problem that helped hypothesize an equilibrium path connecting four pre- and post buckling equilibrium states (solutions). These presented results demonstrate that the proposed procedure can {\color{blue} be adopted to} describe complex post-buckling behaviors of {\color{teal} large structural systems}. Therefore, {\color{teal} this methodology can represent an effective alternative approach to Newton-Raphson strategies for buckling analyses of truss structures or combined with them to improve the robustness and accuracy of structural analyses.


This research provided solutions to this challenging nonlinear, multimodel, unconstrained minimization optimization problem, which the optimization research community can consider as a testbench to test and evaluate new optimization algorithms.} This research was able to solve this multimodel optimization problem by running various instances. However, one would ideally run a single instance to find as many solutions as possible. For example, the population diversity of DE (the best performing algorithm for this set of problems) was extremely low, i.e., almost all individuals in the population produce extremely similar solutions. Hence, only one point on the equilibrium path could be considered in one instance of a run. Therefore, this problem can be presented as a test problem to assess the quality of an algorithm's diverse solutions.

Additionally, the number of iterations required for solving problems with increasing DoF (free variables) was exponentially increasing. For example, space truss structure optimization of 4 variables took on average between 150--500 iterations, 16 variables took on average between \num{5000}--\num{10000} iterations (for some solutions, it took larger than \num{10000} iterations depending on the position of solutions on the equilibrium path), {\color{blue} 21 variable took \num{20000}--\num{25000} iterations,} and 30 variables took about \num{100000}--\num{150000} iterations. Hence, the dimension of this problem is also presenting significant challenges to optimization algorithms. In summary, this test problem presents challenges to optimization algorithms to produce accurate and diverse solutions with high convergence speed.

\section{Conclusions}
\label{sec:conclusion}
{\color{blue}
This work presents a novel analysis methodology based on gradient-free optimization algorithms for the nonlinear structural analysis of space trusses. The proposed methodology formulated the optimization of the global equilibrium of the system as a nonlinear, multimodal, unconstrained, continuous optimization problem. This problem is solved within a new effective multistep analysis procedure providing the nonlinear load path of the system. The proposed methodology can represent an effective method, alternative to Newton-Raphson procedures or combined with them for nonlinear post-buckling analysis of real structural systems.

The search landscape of this problem and the interaction between the free (displacement) variables pose significant challenges to existing continuous optimization algorithms to produce diverse and accurate solutions with high convergence speed. In this research, a number of strategies for search domain decomposition are presented. Consequently, a novel adaptive research space partition method, called hypersphere search algorithm, is proposed and applied to solve this problem. This algorithm iteratively moves through the search landscape of the problems to find as many accurate solutions as possible to provide the complete load-path curve of the structure.

Among the different investigated algorithms, the Differential Evolution (DE) algorithm is identified as the most competitive algorithm to solve this class of problems. Therefore, DE algorithms were applied to solve a test problem concerning a 3D medium-sized reticular beam. As a result, the algorithm optimally produced deformed equilibrium shapes and hypothesis equilibrium path, proving the capability to the procedure to be applied for assessing real structural systems.

}

\small

\bibliographystyle{IEEEtranSA} 
\bibliography{0CivilEng}

\end{document}